\documentclass[12pt]{amsart}
\usepackage{amsmath,amssymb}
\usepackage{amsthm}
\usepackage{epsf}
\usepackage{amssymb}
\usepackage{latexsym}
\usepackage{amsmath}
\newtheorem{theorem}{Theorem} 
\newtheorem{definition}[theorem]{Definition}

\newtheorem{corollary}[theorem]{Corollary}
\newtheorem{lemma}[theorem]{Lemma}
\newtheorem{example}[theorem]{Example}

\newtheorem{remark}[theorem]{Remark}

\newtheorem{notation}[theorem]{Notation}

\def \R{{\mathbb R}}

\def\C{{\mathbb C}}

\def\Z{{\mathbb Z}}
\def\P{{\mathbb P}}

\def\semidirect{\rtimes}

\def \ta{\tau}
\def \ta1{\tau_1}

\def \G{\Gamma}

\def \B{{\mathcal B}}

\def \prodl{\prod\limits}

\def\vK{{van Kampen}}

\newcommand\set[1]{{\{{#1}\}}}

\def\st{{such that }}

\def\isom{{\cong}}

\newcommand\begintable[1][] {{}}

\long\def\forget#1\forgotten{}

\newif\ifXY 
\XYtrue     

\begin{document}

\renewcommand{\subjclassname}{%

    \textup{2000} Mathematics Subject Classification}

\title[Toric varieties - degeneration, fundamental groups]{Degenerations and fundamental groups
related to some special toric varieties}

\author[Meirav~Amram, Shoetsu~Ogata]{Amram Meirav$^1$ and Shoetsu Ogata$^2$}
\stepcounter{footnote}
\footnotetext{This work is partially supported by the Edmund Landau Center for
Research in Mathematical Analysis and Related Areas, sponsored by the
Minerva Foundation (Germany);  DAAD and EU-network HPRN-CT-2009-00099(EAGER);
The Emmy Noether Research Institute for Mathematics and the Minerva
Foundation of Germany, The Israel Science Foundation grant \# 8008/02-3
(Excellency Center ``Group  Theoretic Methods in the Study of Algebraic Varieties").}
\stepcounter{footnote}
\footnotetext{Partially supported by the Ministry of Education, Culture,
Sports, Science and Technology, Japan.}

\address{Meirav Amram, Einstein Institute for Mathematics, Hebrew University, Jerusalem;
Department of Mathematics, Bar-Ilan University, Ramat-Gan, Israel}
\email{ameirav@math.huji.ac.il / meirav@macs.biu.ac.il}
\address{Shoetsu Ogata, Mathematical Institute, Tohoku University, Sendai 980, Japan}
\email{ogata@math.tohoku.ac.jp}


\forget
\begin{abstract}
In this paper we calculate fundamental groups (and some of their quotients)
of complements of four toric varieties branch curves. For these calculations,
we study properties and degenerations
of these toric varieties and the braid monodromies of the branch curves in $\C\P^2$.

The fundamental groups related to the first three toric
varieties turn to be quotients of the Artin braid groups $\mathcal{B}_5$, $\mathcal{B}_6$,
and $\mathcal{B}_4$, while the fourth one is a certain quotient of the group $\tilde{\B}_6 =
\B_6/\langle[X,Y]\rangle$, where $X, Y$ are transversal.

The quotients of all four groups by the normal subgroups generated by the
squares of the standard generators are respectively $S_5, S_6, S_4$ and $S_6$.
We therefore conclude that the fundamental groups of the Galois covers of the four given toric
varieties are all trivial.
\end{abstract}
\forgotten



\date{\today}

\maketitle
\section{Introduction}\label{intro}

    Let $X$ be a projective algebraic surface embedded in a projective space $\C\P^N$.
Take a general linear subspace $V$ in   $\C\P^N$ of dimension $N-3$.  Then the
projection centered at $V$ to $\C\P^2$ defines a finite map
$f: X \to \C\P^2$.  Let $B \subset \C\P^2$ be the branch curve of $f$.
Denote $\pi_1(\C\P^2\setminus B)$ to be {\it the fundamental group of the complement
of the branch curve}.  This group is an invariant of the surface.
Closely related to this group is the affine part $\pi_1(\C^2 \setminus B)$.

In this work we compute the above defined groups, related to four toric varieties.
The first surface is $X_1:=F_1=\P(\mathcal{O}\oplus \mathcal{O}(1))$, the Hirzebruch surface
of degree one in $\C\P^6$ embedded by
the line bundle with the class $s+3g$, where $s$ is the negative
section and $g$ is a general fiber.
The second is $X_2:=F_0=\C\P^1\times \C\P^1$,
the Hirzebruch surface of degree zero in $\C\P^7$ embedded by $\mathcal{O}(1,3)$.
We generalize the results to the case where $X_2$ is embedded in $\C\P^{2n+1}$ by
$\mathcal{O}(1,n)$. The third is $X_3:=F_2=\P(\mathcal{O}\oplus \mathcal{O}(2))$
in $\C\P^5$ embedded by the class $s+3g$.
The fourth is a singular toric surface $X_4$ with one $A_1$ singular point embedded
in $\C\P^6$. $A_1$-singularity is an isolated normal singularity of
dimension two whose resolution consists of one $(-2)$-curve (i.e., a nonsingular rational
curve on a surface with $-2$ as its self-intersection number).
For the first three cases, we use different triangulations of tetragons from
those treated  in \cite{MT4} and \cite{MT5}.

This work fits into the program initiated by Moishezon and Teicher to study
complex surfaces via braid monodromy techniques. They defined the generators
of a braid group from a line arrangement in $\C\P^2$, which is the branch curve of a generic
projection from a union of projective planes \cite{MT4}, namely degeneration.
In order to explain the process of such a degeneration, they used
schematic figures consisting of  triangulations of triangles and tetragons (\cite{MT1},
\cite{ZcZ}, \cite{MT4}). Moishezon and Teicher studied the cases when $X$ is the projective
plane embedded by $\mathcal{O}(3)$ \cite{MT4}, or when $X$ are Hirzebruch surfaces $F_k(a,b)$
for $a, b$ relatively prime \cite{MRT}. Later works compute the group $\pi_1(\C\P^2 \setminus B)$
related to $K3$ surfaces \cite{K3-a}, $\C\P^1 \times T$ where $T$ is a complex torus
(\cite{Gol1}, \cite{Gol2}), $T \times T$ (\cite{AT-TT1}, \cite{AT-TT2}), and
Hirzebruch surface $F_1(2,2)$ \cite{ATV-H}. A very interesting and helpful work concerning
degenerations, braid monodromy and fundamental groups, was written by Auroux-Donaldson-Katzarkov-Yotov
\cite{denis}.

We consult the above works and give a geometric meaning to
these schematic figures from the point of view of toric geometry (\cite{F}, \cite{Od}).
The work is done along the following lines. First we degenerate $X$ into a union $X_0$ of planes.
Then $X_0$ is composed of $n=\deg(X_0)$
planes. $B_0$ is the union of the intersection lines $1, 2, ..., m$ (as depicted in Figures
\ref{numlines}, \ref{1on3}, \ref{trapez}, \ref{home}). The lines are numerated for
future use. It is very complicated to get a presentation of $\pi_1(\C^2\setminus B)$ directly,
therefore we use the regeneration rules from \cite{MT5} to get a braid monodromy factorization
of $B$ from the one of $B_0$. Then we can use the van Kampen Theorem \cite{vK} to get a
finite presentation of $\pi_1(\C^2\setminus B)$ with generators $\Gamma_1, \Gamma_{1'},
\dots, \Gamma_{m}, \Gamma_{m'}$ ($2m$ is the degree of $B$). A presentation of
$\pi_1(\C\P^2\setminus B)$ is obtained by adding the projective relation
$\Gamma_{m'} \Gamma_{m} \cdots \Gamma_{1'} \Gamma_{1}=e$.
The reader might want to check the papers \cite{AT-TT2}, \cite{Gol1}, \cite{Cox} and \cite{ATV-H}
in order to get the feeling of the type of presentations we are dealing with.

Artin \cite{artin} defined the braid group $\mathcal{B}_n$ with $n-1$ generators $\{\sigma_i\}$
and with the following relations
\begin{eqnarray}
{} \sigma_i \sigma_j & =& \sigma_j \sigma_i \ \ \mbox{for $|i-j|>1$}\\
{} \sigma_i \sigma_{i+1} \sigma_i &=& \sigma_{i+1} \sigma_i \sigma_{i+1}.
\end{eqnarray}
The main results in this work, related to $X_1,X_2, X_3$, appear in Theorems \ref{pres1},
\ref{cp1on3} and \ref{pretrapez}:
\begin{itemize}
\item
$\pi_1(\C\P^2 \setminus B_1) \ \isom \ {\mathcal{B}_5}/{\langle\G_{4}^2  \G_{3} \G_{2}
\G_{1}^2 \G_{2} \G_{3}\rangle}$,
\item
$\pi_1(\C\P^2 \setminus B_2) \ \isom \ {\mathcal{B}_6}/{\langle\G_{3} \G_4 \G_{5}^{2} \G_{4} \G_3 \G_{2}
\G_1^2 \G_{2}\rangle}$,
\item
$\pi_1(\C\P^2 \setminus B_3) \ \isom \ \mathcal{B}_4/\langle\G_{2} \G_{3}^2 \G_{2} \G_1^2\rangle$.
\end{itemize}

\begin{remark}
The groups $\pi_1(\C\P^2 \setminus B_i)$ are in fact the braid group of points on the sphere.
A general geometric interpretation is the following.
The surfaces $X_i$ $(i=~1, 2, 3)$ are ruled surfaces, and if $p$ is any point
of $\C\P^2$ outside the branch curve, then its $N$ preimages in $X_i$ $(N=5, 6, 4)$
project to distinct points of $\C\P^1$; this gives a homomorphism from $\pi_1(\C\P^2 \setminus {B_i})$
to $B_N(\C\P^1)$.
\end{remark}

The result related to $X_4$ appears in Theorem \ref{prehome}:
\begin{itemize}
\item
{\em $\pi_1(\C\P^2 \setminus B_4)$ is isomorphic to a quotient of the group $\tilde{\B}_6 =
\B_6/\langle[X,Y]\rangle$ $(X, Y$ are transversal$)$ by $\langle(\ref{anilo})\rangle$}.
\end{itemize}

\smallskip
In this work we are also interested in two important quotient groups.
The first one  $\Pi_{(B)}=\pi_1(\C\P^2\setminus B)/\langle\Gamma_{i}^2,
\Gamma_{i'}^2\rangle$ is defined to be a quotient of $\pi_1(\C\P^2\setminus B)$ by the normal subgroup
generated by the squares of the generators.
This group is a key ingredient in studying invariants of $X$, and in particular
$\pi_1(\C\P^2\setminus B)$.
The braid monodromy technique of Moishezon-Teicher enables one to compute $\pi_1(X_{Gal})$,
the fundamental group of a Galois cover $X_{Gal}$ of $X$, from $\Pi_{(B)}$.
In particular, they showed that there is a natural map from $\Pi_{(B)}$
to the symmetric group $S_n$, where $n$ is the degree of $X$,
and $\pi_1(X_{Gal})$ is the kernel of this homomorphism.
Moishezon-Teicher proved in \cite{ZcZ} that for $X= \C\P^1 \times \C\P^1$ the group
$\pi_1(X_{Gal})$ is a finite abelian group on $n-2$ generators, each of
order g.c.d.$(a,b)$ ($a$ and $b$ are the parameters of the embedding).
In \cite{Gol1} the treated surface is $X=\C\P^1 \times T$ ($T$ is a complex torus) and
$\pi_1(X_{Gal})=\Z^{10}$. In \cite{Gol2} the same surface was embedded in $\C\P^{2n-1}$ and
$\pi_1(X_{Gal})=\Z^{4n-2}$. In \cite{Cox} and \cite{ATV} the surface $X=T \times T$
is studied, and $\pi_1(X_{Gal})$ is nilpotent of class $3$. In \cite{ATV-H}
this group was computed for the Hirzebruch surface $F_1(2,2)$ and it is $\Z_2^{10}$.

It turns out in this paper (Theorems \ref{pres1}, \ref{cp1on3}, \ref{pretrapez},
\ref{prehome}) that
\begin{itemize}
\item
{\em The group $\Pi_{(B_i)}$ is isomorphic to  $S_5, S_6, S_4, S_6$ for $i=1, 2, 3, 4$, respectively.}
\end{itemize}
Hence we have
\begin{corollary}
The fundamental group $\pi_1({(X_i)}_{Gal})$ is trivial for $i=1, 2, 3, 4$.
\end{corollary}

The second group is a Coxeter group  $C = \Pi_{(B)}/\langle\G_i=\G_{i'}\rangle$, defined as a quotient of
$\Pi_{(B)}$ under identification of pairs of generators, see \cite{RTV}. It is still unclear whether $C$,
introduced here, is an invariant of the surface or of the branch curve.
It might be conjectured that there exists a dependence on the choice of a pairing between geometric
generators $\G_j$ and $\G_{j'}$ (and hence on the choice of a degeneration to a union of planes).
It turns out that $C$ is isomorphic to a symmetric group $S_n$ for
Hirzebruch surfaces (\cite{ATV-H}, \cite{MRT}) and $\C\P^1 \times \C\P^1$ (\cite{MT1}, \cite{ZcZ}).
The cases of $\C\P^1 \times T$ (\cite{Gol1}) and $T \times T$ (\cite{Cox}) are the first examples
in which $C$ is a larger group, namely  $C \ \isom \ \Z_5 \semidirect S_{6}$ and $C \ \isom \ K_C \semidirect S_{18}$
($K_C$ is a central extension of $\Z^{34}$ by $\Z$), respectively.

Then we have
\begin{corollary}
The group $C_i$ is isomorphic to $S_5, S_6, S_4, S_6$ for $i=1, 2, 3, 4$, respectively.
\end{corollary}

The paper is divided as follows. In Section \ref{sec2} we study
degeneration of toric varieties.  In Section \ref{sur1} we compute the
requested groups related to the toric varieties $X_1, X_2$ and $X_3$,
and in Section \ref{x4} we compute the ones related to $X_4$.

\medskip

{\bf Acknowledgements.}
This research was initiated while the two authors were at the Mathematics
Institute, Erlangen - N\"urnberg University, Germany. They both wish to thank the
Institute for its hospitality. The first author wishes to thank her hosts Profs.
W. Barth and H. Lange. The second author would like to thank his host Prof. H. Lange.
The first author also would like to thank the Department of Mathematics,
Bar-Ilan University and her host Prof. Mina Teicher and the Einstein Institute
for Mathematics, Jerusalem, and her host Prof. Ruth Lawrence-Neumark.

Both authors wish to thank the referee for very important remarks, Denis Auroux for fruitful
discussions and guidance in preparing the bibliography, and to Michael Friedman for helpful comments.

\section{Degeneration of toric surfaces}\label{sec2}
In their process to calculate the braid monodromy, Moishezon and Teicher
studied the projective degeneration of $V_3=(\C\P^2, {\mathcal O}(3))$ \cite{MT4}
and Hirzebruch surfaces \cite{MRT}.
Since $\C\P^2$ and the Hirzebruch surfaces are toric surfaces, we shall
describe the projective degeneration of toric surfaces in this section.

\subsection{Basic notions}\label{sec2a}
We outline definitions needed in toric geometry and refer to \cite{F} and
\cite{Od} for further statements and proofs.

\begin{definition}{\bf Toric variety.}
A toric variety is a normal algebraic variety $X$ that contains an algebraic torus
$T=({\C}^*)^n$ as a dense open subset, together with an algebraic action
$T\times X \to X$ of $T$ on $X$, that is an extension of the natural action of $T$
on itself.
\end{definition}

Let $M$ be a free $\Z$-module of rank $n$ $(n \ge 1)$ and $M_{\R}:=M \otimes_{\Z} \R$
the extension of the coefficients to the real numbers. Let $T:=\mbox{Spec}{\C}
[M]$ be an algebraic torus of dimension $n$. Then $M$ is considered as the character
group of $T$, i.e.,  $M=\mbox{Hom}_{\mbox{gr}}(T, {\C}^*)$. We denote an element
$m \in M$ by $e(m)$ as a function on $T$, which is also a rational function on $X$.
Let $L$ be an ample line bundle on $X$.  Then we have
\begin{equation}
H^0(X, L) \cong \bigoplus_{m \in P \cap M} {\C}e(m),
\end{equation}
where $P$ is an integral convex polytope in $M_{\R}$ defined as the convex hull \linebreak
Conv$\{m_0, m_1, \dots, m_r\}$ of a finite subset $\{m_0, m_1, \dots, m_r\}
\subset M$. Conversely we can construct a pair $(X, L)$ of a polarized toric
variety from an integral convex polytope $P$ so that the above isomorphism holds
(see \cite[Section 3.5]{F} or \cite[Section 2.4]{Od}).
If an affine automorphism $\varphi$ of $M$ transforms $P$ to $P_1$, then
$\varphi$ induces an isomorphism of polarized toric varieties $(X, L)$ to
$(X_1, L_1)$, where $(X_1, L_1)$ corresponds to $P_1$.

\begin{example}
Let $M={\Z}^2$.  Then $V_3=({\C\P}^2, {\mathcal O}(3))$ corresponds to the
integral convex polytope $P_3:=$\emph{Conv}$\{(0,0), (3,0), (0,3)\}$.
\end{example}
\begin{example}\label{ex:1}
The Hirzebruch surface
$F_d={\P}({\mathcal O}_{\C\P^1} \oplus {\mathcal O}_{\C\P^1}(d))$ of degree $d$
has generators $s, g$ in the Picard group consisting of the negative section $s^2=-d$
and general fiber $g^2=0$. A line bundle $L$ with $[L]=as + bg$ in $\mbox{Pic}
({F}_d)$ is ample if $a>0$, $b>ad$.  Then this pair $({F}_d, L)$
corresponds to $P_{d(a,b)}:=$\emph{Conv}$\{(0,0), (b-~ad,0), (b, a), (0,a)\}$.
\end{example}

Next we consider degenerations of toric surfaces defined by Moishezon-Teicher.
We recall the definition from \cite{MT4}.
\begin{definition}\label{df1}{\bf  Projective degeneration.}
A degeneration of $X$ is a proper surjective morphism with connected fibers
$\pi  : V \rightarrow \C$ from an algebraic variety $V$
 such that the restriction $\pi:
V \setminus \pi^{-1}(0) \rightarrow \C \setminus \{0\}$
is smooth and that $\pi^{-1}(t) \cong X$ for $t\not=0$.

When $X$ is projective with an embedding $k: X \hookrightarrow \C \P ^n$,
 a degeneration of $X$ $\pi  : V \rightarrow \C$ is called {\bf  a projective degeneration}
of $k$ if there exists a  morphism $F  : \ V \rightarrow
\C \P^n \times \C$ \st the restriction
$\ F_t = F \mid_{\pi^{-1}(t)} \ : \pi^{-1}(t) \rightarrow \C \P^n
\times t$ is an embedding of $\pi^{-1}(t)$ for all $t\in \C$
and that $F_1 = k$ under the identification of
$\pi^{-1}(1)$ with $X$.
\end{definition}

Moishezon and Teicher used the
triangulation of $P_3$ consisting of nine standard triangles as a schematic figure
of a union of nine projective planes \cite{MT4}. In the theory of toric varieties,
however, the lattice points $P_3 \cap M$ correspond to rational functions of degree $3$
on $V_3\cong \C\P^2$.  Let $m_0=(0,0), m_1=(1,0), m_2=(0,1), \dots, m_9=(0,3)\in \Z^2$.
Then we may write $e(m_0)=x_0^3, e(m_1)=x_0^2x_1, e(m_2)=x_0^2x_2, \dots, e(m_9)=x_2^3$
with a suitable choice of the homogeneous coordinates of $\C\P^2$.
The Veronese embedding $V_3 \hookrightarrow \C\P^9$ is given by $z_i =e(m_i)$
for $i=0,1, \dots, 9$ with the homogeneous coordinates $[z_0: z_1: \dots :z_9]$
of $\C\P^9$.  Let $P_1:= \mbox{Conv}\{(0,0), (1,0), (0,1)\}$, which corresponds
to
$(\C\P^2, {\mathcal O}(1))$.  The subset $P_1 \subset P_3$ corresponds to the
linear subspace $\{z_3=\cdots =z_9 =0\} \subset \C\P^9$.
Thus a  triangulation of $P_3$ into a union of nine standard triangles means
the subvariety of dimension two consisting of the union of nine projective planes
in $\C\P^9$ and each standard triangle defines a linear subspace of dimension two
with corresponding coordinates.

\subsection{Constructing the degeneration of toric surfaces}\label{sec2b}
In the following
we construct a semistable degeneration of toric surfaces according to Hu \cite{Hu}.
Let $M=\Z^2$. Let $P$ be a convex polyhedron in $M_{\R}$ corresponding to a
polarized toric surface $(X, L)$.  The lattice points $P \cap M$ define the embedding
$\varphi_L: X\to \P(\Gamma(X, L))$. Let $\Gamma$ be a triangulation of $P$ consisting
of standard triangles with vertices in $P \cap M$.  Let $h: P\cap M \to \Z_{>0}$ be a
function on the lattice points in $P$ with values in positive integers.
Let $\tilde M=M \oplus \Z$ and let $\tilde P=\mbox{Conv}\{(x,0), (x, h(x));
x \in P \cap M\}$ the integral convex polytope in $\tilde M_{\R}$.
We want to choose $h$ to satisfy the conditions that $(x, h(x))$ for $x\in P\cap M$
are vertices of $\tilde P$ and that for each edge in $\Gamma$ joining $x$ and $y\in P\cap M$
there is an edge joining $(x, h(x))$ and $(y, h(y))$ as a face of $\partial \tilde P$.
We say that $\tilde P$ realizes the triangulation $\Gamma$ if these conditions are
satisfied.
Now we assume that
$\tilde P$ realizes the triangulation $\Gamma$.  Then $\tilde P$ defines a polarized
toric 3-fold $(\tilde X, \tilde L)$. From the construction, $\tilde X$ has a fibration
$p:\tilde X\to \C\P^1$ satisfying that $p^{-1}(t)\cong X$ with $t \not = 0$ and that
$p^{-1}(0)$ is a union of projective planes.
Furthermore we see that $p^{-1}({\C\P}^1 \setminus \{0\}) \cong {\C}\times X$.
Thus the flat family $p:\tilde X\to \C\P^1$ gives a degeneration of $X$ into a union
of projective planes with the configuration diagram $\Gamma$.
Hu treats only nonsingular toric varieties of any dimension.
The difficulty of this construction is to find a triangulation $\Gamma$.
Here we restrict ourselves to toric surfaces.  Then we can find a
triangulation for any integral convex polygon $P$.

\begin{example}
Let $m_0=(0,0), m_1=(1,0), m_2=(0,1), m_3=(1,1) \in M=\Z^2$.
Let $P=$\emph{Conv}$\{m_0, m_1, m_2, m_3\}$.  Then $P$ defines
the polarized surface $(X=\C\P^1 \times \C\P^1, {\mathcal O}(1,1))$.
Let $\Gamma$ be the triangulation of $P$ defined by adding the edge
connecting $m_1$ and $m_2$.  Define $h(m_0)=h(m_3)=1, h(m_1)=h(m_2)=2$.
Let $\tilde M:= M \oplus \Z$.
Set $m_i =(m_i,0)$ and $m_i^+ =(m_i, h(m_i))$ for $i=0, \dots, 3$
and $m_4=(1,0,1), m_5=(0,1,1)$ in $\tilde{M}$.
Then the integral convex polytope $\tilde{P} := \mbox{Conv}\{
m_0, \dots, m_3, m_0^+, \dots, m_3^+\}$ in $\tilde{M}$ defines
the polarized toric 3-fold $(\tilde{X}, \tilde{L})$.
By definition, $\tilde{X}$ has a fibration $p:~\tilde{X} \to~\mathbb{C}\mathbb{P}^1$.
The global sections of $\tilde{L}$ defines an embedding of $\tilde{X}$
as follows:  Let $[z_0: \dots: z_9]$ be the homogeneous coordinates
of $\mathbb{C}\mathbb{P}^9$.  The equations $z_i =e(m_i)$ for $i=0, \dots, 5$
and $z_{6+j}=~e(m_j^+)$ for $j=0, \dots, 3$  define the embedding
$\tilde{X} \to \mathbb{C}\mathbb{P}^9$.
The fiber $p^{-1}(\infty)$ is given by $\{z_0z_3=z_1z_2, z_4=\dots =z_9=0\}$
which is isomorphic to $X  \subset~\mathbb{P}(\Gamma(X, \mathcal{O}(1,1)))
\cong \mathbb{C}\mathbb{P}^3 =\{z_4=\dots =z_9=0\}$,
and the fiber $p^{-1}(0)$ is given by
$\{z_6z_9=0, z_0=\dots =z_5=0\}$ which is a union of two projective planes
in $\mathbb{C}\mathbb{P}^3 \cong\{z_0=\dots =z_5=0\}$.
\end{example}

\begin{lemma}
The line bundle $\tilde L$ on $\tilde X$ is very ample.
\end{lemma}

\begin{proof}
Let $m_1, m_2,  m_3 \in P \cap M$ be three vertices of a standard triangle in the
triangulation $\Gamma$ of $P$. Set $m_i^-=(m_i, 0), m_i^+=(m_i,h(m_i))$ in
$\tilde M_{\R}$ for $i=1,2,3$. Denote $Q=\mbox{Conv}\{m_i^{\pm}; i=1,2,3\}$ the
integral convex polytope with vertices $\{m_i^{\pm}; i=1,2,3\}$. Then we divide
$\tilde P$ into a union of triangular prisms like $Q$. We can divide $Q$ into a union
of standard 3-simplices.
We may assume $h(m_1)\ge h(m_2)\ge h(m_3)$
by renumbering $m_i$ if necessary. Then we can
divide $Q$ into a union of $Q_0=\mbox{Conv}\{m_1^+, m_2^+, m_3^+, (m_1, h(m_1)-1)\}$
and $Q_1=\mbox{Conv}\{m_1^-, (m_1, h(m_1)-1), m_2^{\pm}, m_3^{\pm}\}$.
Here $Q_0$ is a standard 3-simplex and $Q_1$ has a similar shape to $Q$ but less
volume than that of $Q$. Thus we obtain a division of $\tilde P$ into a union of
standard 3-simplices. This is not always triangulation of $\tilde P$, but this gives a
covering of $\tilde P$ consisting of standard 3-simplices.  From the theory of
polytopal semigroup ring (see, for instance, \cite{BGT} and \cite{Str}), we see that
$\tilde L$ is simply generated, hence very ample.
\end{proof}

We claim that $\tilde X$ also defines a projective degeneration of $(X, L)$.
Denote \linebreak $\Phi:=\varphi_{\tilde{L}}: \tilde{X} \longrightarrow \mathbb{P}
(\Gamma(\tilde{X}, \tilde{L})) =: \mathbb{P}$ the morphism defined by global
sections of $\tilde{X}$. We see that $p^{-1}(t)\cong X$ for $t\not=0$ with
$[1:t] \in \mathbb{C}\mathbb{P}^1$ and that $p^{-1}(\infty)\cong X$ and $p^{-1}(0)$
are $T$-invariant reduced divisors. Thus the restriction maps $\Gamma(\tilde{X}, \tilde{L})
\longrightarrow \Gamma(p^{-1}(\infty), \tilde{L}|_{p^{-1}(\infty)})\cong \Gamma(X, L)$
and $\Gamma(\tilde{X}, \tilde{L}) \longrightarrow \Gamma(p^{-1}(0), \tilde{L}|_{p^{-1}(0)})$
are surjective.  From the construction of $\tilde{P}$, we see that $\dim\Gamma(X,L)
 =\dim \Gamma(p^{-1}(0), \tilde{L}|_{p^{-1}(0)})$.
Since $p^{-1}(\mathbb{C}\mathbb{P}^1 \setminus \{0\}) \cong X\times \mathbb{C}$,
we have $\tilde{L}|_{p^{-1}(t)}\cong L$ for $t\not=0$. Then $F:= \Phi \times p:
\tilde{X} \longrightarrow \mathbb{P}\times \mathbb{C}\mathbb{P}^1$ is a projective
degeneration of $k: X \longrightarrow \mathbb{P}(\Gamma(X,L)) \hookrightarrow \mathbb{P}$.

\forget
We claim that $\tilde X$ also defines a projective degeneration of $(X, L)$.
Denote $\varphi_{\tilde L}: \tilde X \to {\P}(\Gamma(\tilde X, \tilde L))$ the
embedding defined by global section of $\tilde L$.  Let $P \cap M=\{m_0, m_1, \dots,
m_r\}$. Set $m_i^-=(m_i, 0), m_i^+=(m_i, h(m_i))$ in $\tilde M$ for $i=0, \dots, r$.
Now we define a morphism $\Phi: \tilde X \to \C\P^r$ as the $i$-th homogeneous
coordinate of $\Phi(x)$ for $x \in \tilde X$ is $\lambda e(m_i^+) +\mu e(m_i^-)$,
where $e(m_i^{\pm})$ are homogeneous coordinates of $\varphi_{\tilde L}(x)$ in
${\P}(\Gamma(X, L))$ and $[\lambda: \mu]$ is the homogeneous coordinate of $p(x)$
in ${\C\P}^1$. Then the restriction of $\Phi$ to $p^{-1}(\lambda=0)$
coincides with $\varphi_L$ and $\Phi(p^{-1}(\mu=0))$  is a union of projective planes
corresponding to the triangulation $\Gamma$ of $P$.
Let $V:= \tilde X\setminus p^{-1}(\lambda=0)$.  Then $F:= \Phi\times p: V \rightarrow
\C\P^r\times \C$ is a projective degeneration of $\varphi_L$.
\forgotten

\begin{theorem}
Let $P$ be an integral convex polyhedron of dimension $2$ corresponding to a polarized
toric surface $(X, L)$ and let $\Gamma$ a triangulation of $P$ consisting of standard
triangles with vertices in $M$. Assume that $\tilde P$ is an integral convex polytope
in $\tilde M_{\R}$ realizing the triangulation $\Gamma$. Then $\tilde P$ defines
a polarized toric 3-fold $(\tilde X, \tilde L)$, which gives a projective degeneration
of $(X, L)$ to a union of projective planes.
\end{theorem}

\subsection{Degeneration of the four toric surfaces}\label{sec2c}
In this paper we study four degenerations of polarized toric surfaces, each one of which
is defined by integral convex polygon $P$.
We choose a triangulation $\Gamma$ for each $P$ and define a function $h:P\cap M \rightarrow
\Z_{\ge0}$ so that the integral convex
polytope $\tilde P$ of dimension $3$ should realize the triangulation $\Gamma$ of $P$.

The first surface is the Hirzebruch surface $X_1:=F_1$ of degree one embedded in
$\C\P^6$ by the very ample line bundle $L_1$ whose class is $s+3g$, where
$s$ is the negative section and $g$ is a general fiber.
We mentioned this surface as a polarized
toric surface in Example~\ref{ex:1}, which corresponds to the integral convex polygon
$P_{1(1,3)}$ in $M=\Z^2$.
Let $m_i=(i,0)$ for $i=0,1,2,3$ and $m_j=(j-3,1)$ for
$j=4,5,6$.  Then $P_{1(1,3)} = \mbox{Conv}\{m_0, m_3, m_4, m_6\}$.
Let $\Gamma_1$ be the triangulation of $P_{1(1,3)}$ obtained by adding the edges
$\bar{m_1m_4}, \bar{m_2m_4}, \bar{m_2m_5}, \bar{m_3m_5}$, see Figure \ref{numlines}.  This triangulation
is slightly different from the one treated in \cite{MT4}.  We define a function
$h_1: P_{1(1,3)} \cap M \rightarrow \Z_{>0}$ as $h_1(m_0)= h_1(m_6)= 1,
h_1(m_1)=h_1(m_3)= h_1(m_4)=h_1(m_5)=3, h_1(m_2)=4$. Then we can define an integral
convex polytope $\tilde P$ in $\tilde M = M\oplus \Z$ realizing the triangulation $\Gamma_1$
of $P_{1(1,3)}$.  Hence we have a projective degeneration of
$\varphi_1:=\varphi_{L_1}: F_1 \hookrightarrow \C\P^6$.

The second surface is $X_2:=\C\P^1\times \C\P^1$ embedded in $\C\P^7$ by ${\mathcal O}(3,1)$.
This embedded toric surface corresponds to the convex polygon $P_{3,1}:=\mbox{Conv}\{(0,0),
(3,0), (0,1), (3,1)\}$ in $M=\Z^2$. Let $m_i=(i,0)$ for $i=0,1,2,3$ and $m_j=(j-4, 1)$
for $j=4,5,6,7$.   Then $P_{3,1}=\mbox{Conv}\{m_0, m_3, m_4, m_7\}$.
Let $\Gamma_2$ be the triangulation of $P_{3,1}$ obtained by adding the edges
$\bar{m_0m_5}, \bar{m_1m_5}, \bar{m_1m_6}, \bar{m_2m_6}, \bar{m_2m_7}$, see Figure \ref{1on3}.
We define a function $h_2: P_{3,1}\cap M \rightarrow \Z_{>0}$ as $h_2(m_4)=1,
h_2(m_0)=h_2(m_3)=3, h_2(m_5)=h_2(m_7)=4, h_2(m_1)=h_2(m_2)= h_2(m_6)=5$.
Then we have a projective degeneration of $\varphi_2:= \varphi_{\mathcal O(3,1)}:
\C\P^1\times \C\P^1 \hookrightarrow \C\P^7$ corresponding to the triangulation
$\Gamma_2$.

The third surface is the Hirzebruch surface $X_3:= F_2$ of degree $2$ embedded in
$\C\P^5$ by the ample line bundle $L_2$ whose class is $s+3g$.  The corresponding
polygon is $P_{2(1,3)}$.  Let $m_i=(i,0)$ for $i=0,1,2,3$ and $m_j=(j-3,1)$ for
$j=4,5$ in $M=\Z^2$.  Then $P_{2(1,3)} = \mbox{Conv}\{m_0, m_3, m_4, m_5\}$
up to affine automorphism of $M$.  Let $\Gamma_3$ be the triangulation of $P_{2(1,3)}$
obtained by adding the edges $\bar{m_1m_4}, \bar{m_1m_5}, \bar{m_2m_5}$, see Figure~\ref{trapez}.
We define a function $h_3: P_{2(1,3)}\cap M \rightarrow \Z_{>0}$ as
$h_3(m_0)=h_3(m_3)=1, h_3(m_4)=3, \linebreak h_3(m_1)=h_3(m_2)=h_3(m_5)=4$.
Then we have a projective degeneration of \linebreak $\varphi_3:=\varphi_{L_2}:
F_2 \hookrightarrow \C\P^5$ corresponding to the triangulation $\Gamma_3$.

The last surface is a singular toric surface $X_4$ embedded in $\C\P^6$ corresponding to
the polygon $P_4:=\mbox{Conv}\{(0,0), (2,0), (0,1), (1,2), (2,1)\}$.
Let $m_i=(i,0)$ for $i=0,1,2$, $m_j=(j-3,1)$ for $j=3,4,5$ and  $m_6=(1,2)$.
Let $\Gamma_4$ be the triangulation of $P_4$ obtained by adding the edges
$\{\bar{m_i m_4}, \bar{m_1 m_j} ; i=1,3,5,6 \ \mbox{and} \ j=3,5\}$, see Figure~\ref{home}.
We define a function $h_4: P_4 \cap M \rightarrow \Z_{>0}$ as
$h_4(m_0)=h_4(m_2)=1, h_4(m_1)=h_4(m_3)=h_4(m_5)=h_4(m_6)=3, h_4(m_4)=4$.
Then we have a projective degeneration of $\varphi_4: X_4 \hookrightarrow
\C\P^6$ corresponding to the triangulation $\Gamma_4$.

\section{The surfaces $X_1$, $X_2$ and $X_3$}\label{sur1}
In this section we compute the groups $\pi_1(\C^2 \setminus B_i)$,
$\pi_1(\C\P^2 \setminus B_i)$, and  $\Pi_{(B_i)}$ for $i=1, 2, 3$.
Zariski \cite{Z} investigated indirectly complements of the types of curves as
$B_1, B_2$ and $B_3$. We compare our methods and results to those of Zariski.

Using degenerations of toric varieties, such as those that we have here, makes these special
cases of a more general theory, rather than isolated examples.
Having the degenerations of $X_1$, $X_2$ and $X_3$, we project them onto
$\C\P^2$ and get line arrangements. By the Regeneration Lemmas of
Moishezon-Teicher \cite{MT3}, the diagonal lines regenerate to conics which are tangent
to the lines with which they intersect. When the rest of the lines regenerate, each tangency
(the point of tangency of line and conic)
regenerates to three cusps. We end up with cuspidal curves $B_i$, $i=1, 2, 3$.
The existence of nodes in these curves depends on the existence of the `parasitic
intersections' (projecting the degenerations onto $\C\P^2$ causes extra intersections).
By the braid monodromy techniques and regeneration rules of Moishezon-Teicher
(\cite{MT3}, \cite{MT5}), we get the related braid monodromy factorizations
(by \cite{MT2}, each braid of a parasitic intersection, say $Z_{ij}^2$, regenerates to
$Z^2_{ii', jj'}$ in the factorizations), see Notation \ref{nodes}.
We do not use properties of braid groups, but rather the definition of the factorization \cite{MT2},
from which the van Kampen Theorem \cite{vK} for cuspidal curves gives a complete set of
relations for the fundamental groups $\pi_1(\C^2 \setminus B_i)$.

Zariski \cite{Z} gets a collection of local relations without using degeneration
and regeneration, as follows. He uses properties of curves to conclude relations
for certain groups, called the Poincar\'e groups (contemporary fundamental groups).
He defines the class of Poincar\'e groups $G_n$, which practically coincides with the
Artin braid groups \cite{artin}.
A group of type $G_n$ is also a group of automorphism classes of a sphere with $n$ holes
(the points $P_1, \dots, P_n$ are removed), see \cite{mag}.
For generators $g_1, \dots, g_{n-1}$ ($g_1$ connects $P_1$ and $P_2$,
$g_2$ connects $P_2$ and $P_3$, etc.),  Zariski proves that
\begin{equation}\label{comG}
g_i g_j = g_j g_i \ \ \mbox{$\vert i-j \vert \neq 1$}
\end{equation}
\begin{equation}\label{tripG}
g_i g_{i+1} g_i = g_{i+1} g_i g_{i+1}
\end{equation}
\begin{equation}
g_1 g_2 \cdots g_{n-2} g^2_{n-1}  g_{n-2} \cdots g_2 g_1=e
\end{equation}
constitute a complete set of generating relations of $G_n$.
He denotes a rational curve with degree $n$ and $k$ cusps as $(n,k)$.
He shows how the individual generating relations of $G_n$ correspond to
the singularities of a maximal cuspidal curve $(2n-2,3(n-2))$ with
$2(n-2)(n-3)$ nodes. The $(n-2)(n-3)/2$ commutativity relations
(\ref{comG}) are the typical relations at nodes, while the $n-2$
relations (\ref{tripG}) are the typical cusp relations \cite{Z2}.

Concerning the results, the cuspidal curves $B_1, B_2$ and $B_3$
($(8,9), (10,12)$ and $(6,6)$, respectively) fulfill the above statements and
they are maximal. Therefore, Zariski gets the groups  $G_5, G_6$ and $G_4$, respectively.
Here the results related to $X_1, X_2$ and $X_3$, turn out to be the ones
of Zariski, i.e., $\pi_1(\C\P^2 \setminus B_i)$ is a braid group of points on a sphere.

Since we use the degeneration on toric varieties, which is different from that which Zariski did,
it would be worth to give a proof for the groups related to $X_1$. The ones related to
$X_2$ and $X_3$ are computed in a similar way, and therefore the proofs are omitted.

\begin{remark}
A braid monodromy factorization $\Delta^2$ should be written as a product of factors
in an actual order $($see \cite{MT5}$)$. Since our goal is to compute fundamental groups,
the order of the factors does not matter. Here we  list the monodromies with an unmeaningful order,
and concentrate on finding the relations in the groups by applying the van Kampen Theorem on the monodromies.
\end{remark}

\subsection{The surface $X_1$}
Let $X_1 = F_1(3,1)$ be the Hirzebruch surface, as defined in Section \ref{sec2}.
The construction of the degeneration of Hirzebruch surfaces of type $F_1(p,q)$
(for $p>q \geq 2$) appears in \cite{MII} and \cite{M2}. In \cite{denis},
Auroux-Donaldson-Katzarkov-Yotov dedicate Section 6.2 to the construction of
degeneration of $F_1$ surfaces and to the presentations of the fundamental groups of
complements of branch curves.

\begin{figure}[h]
\begin{minipage}{\textwidth}
\begin{center}
\epsfbox{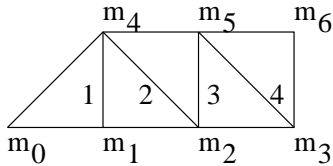} \caption{Degeneration of
$X_1$}\label{numlines}
\end{center}
\end{minipage}
\end{figure}

The degeneration of $X_1$ into a union of five planes ${(X_1)}_0$ is embedded in $\C\P^6$.
The numeration of lines is fixed according to the numeration of the vertices
in Section \ref{sec2}, see Figure \ref{numlines}.
Note that each of the points $m_2, m_4, m_5$ is contained in three distinct planes,
while each of $m_1, m_3$ is contained in two planes.

Take a generic projection $f_1: X_1 \to \C\P^2$.
The union of the intersection lines is the ramification locus $R_0$ in ${(X_1)}_0$ of
$f_1^0 \colon {(X_1)}_0 \to \C\P^2$. Let ${(B_1)}_0 = f_1^0(R_0)$
be the degenerated branch curve. It is a line arrangement,
${(B_1)}_0=\bigcup \limits^{4}_{j=1}L_j$.

Denote the singularities of ${(B_1)}_0$ as $f_1^0(m_i)=m_i, i=1, \dots, 5$
(the points $m_0, m_6$ do not lie on numerated lines, hence they are not
singularities of ${(B_1)}_0$). The points $m_1$ and $m_3$ (resp. $m_2, m_4, m_5$) are called
$1$-points (resp. $2$-points). They were studied in \cite{Gol1}, \cite{ATV-H}, \cite{MT1} and
\cite{MT5}. Other singularities may be the parasitic intersections.

\begin{figure}
\begin{minipage}{\textwidth}
\begin{center}
\epsfbox{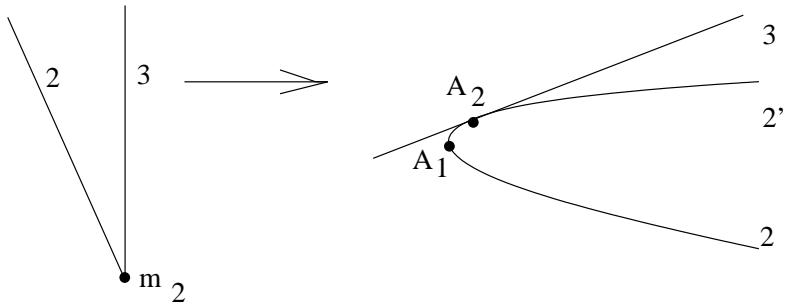} \caption{Regeneration around the point
$m_2$}\label{phi3}
\end{center}
\end{minipage}
\end{figure}

The regeneration of ${(X_1)}_0$ induces a regeneration of ${(B_1)}_0$ in such
a way that each point on the typical fiber, say $c$, is replaced by two close
points $c, c'$. The regeneration occurs as follows.
We regenerate in a neighborhood of $m_1, m_3$ to get conics.
Now, in a neighborhood of  $m_2, m_4, m_5$, the diagonal line regenerates to a
conic \cite[Regenerations Lemmas]{MT3}), which is tangent to the line it intersected
with \linebreak\cite[Lemma 1]{MT5}. See Figure \ref{phi3} for the regeneration around $m_2$.
When the line regenerates, the tangency regenerates into three cusps,
\cite[Regeneration Lemmas]{MT3}.

The resulting curve $B_1$ has degree $8$ and nine cusps.
The intersection points of the curve
with a typical fiber are $\{1, 1', \dots, 4, 4'\}$.
We are interested in the braid monodromy factorization of $B_1$,
the groups $\pi_1(\C^2 \setminus B_1), \pi_1(\C\P^2 \setminus B_1)$ and $\Pi_{(B_1)}$.

\begin{notation}\label{nodes}
We denote by $Z_{i \; j}$ the counterclockwise half-twist of $i$ and $j$
along a path below the real axis. Denote by $Z^2_{i,j \; j'}$ the product
$Z_{i \; j'}^2 \cdot Z_{i \; j}^2$, and by $Z^2_{i \; i', j \; j'}$ the product
$Z_{i', j \; j'}^2 \cdot Z_{i, j \; j'}^2$.
Likewise, $Z^{3}_{i,j \; j'}$ denotes the product $(Z^3_{i \; j})^{Z_{j \; j'}} \cdot (Z^3_{i \; j})
\cdot (Z^3_{i \; j})^{{Z_{j \; j'}}^{-1}}$. Conjugation of braids is defined as $a^b = b^{-1}ab$.
\end{notation}

\begin{theorem}\label{del8}
The braid monodromy factorization of the curve $B_1$ is the product of
\begin{eqnarray}
{}\varphi_{m_1} & = & Z_{1 \; 1'} \label{m1}\\
{}\varphi_{m_2} & = & Z^3_{2', 3 \; 3'} \cdot {Z_{2\; 2'}}^{Z^2_{2', 3 \; 3'}} \label{m2}\\
{}\varphi_{m_3} & = & Z_{4 \; 4'} \label{m3}\\
{}\varphi_{m_4} & = & Z^3_{1 \; 1', 2} \cdot {Z_{2\; 2' }}^{Z^2_{1 \; 1', 2}} \label{m4}\\
{}\varphi_{m_5} & = & Z^3_{3\;3', 4} \cdot {Z_{4\; 4' }}^{Z^2_{3\; 3', 4}} \label{m5}
\end{eqnarray}
and the parasitic intersections braids
\begin{equation}
{Z^2}_{\hspace{-.1cm}{1 \; 1',3 \; 3'}}, \; {Z^2}_{1\; 1',4 \; 4'}, \;
{Z^2}_{2\; 2',4 \; 4'}.\label{par}
\end{equation}
\end{theorem}

\begin{proof}
The monodromies (\ref{m1}) and (\ref{m3}) are derived from the regenerations around \linebreak $1$-points,
and the ones related to $2$-points are (\ref{m2}), (\ref{m4}), (\ref{m5}),
see for example the braids of $\varphi_{m_2}$ in Figure \ref{3phi}.
\begin{figure}[h]
\begin{minipage}{\textwidth}
\begin{center}
\epsfbox{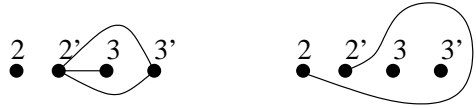} \caption{The braids of
$\varphi_{m_2}$}\label{3phi}
\end{center}
\end{minipage}
\end{figure}
The parasitic intersections were formulated in \cite{MT2}, these are the intersections
of the lines $L_1$ and $L_3$, $L_1$ and $L_4$, $L_2$ and $L_4$.
See Figure \ref{ci}.
\begin{figure}[h!]
\begin{minipage}{\textwidth}
\begin{center}
\epsfbox{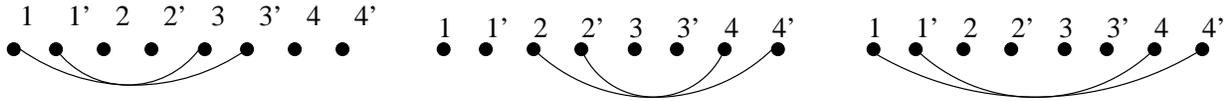} \caption{Parasitic intersections braids in the
factorization of $B_1$}\label{ci}
\end{center}
\end{minipage}
\end{figure}

Summing the degrees of the braids gives $56$. Since the degree of the factorization
is $56$ \cite[Cor. V.2.3]{MT2}, no other braids are involved.
\end{proof}

\smallskip

\begin{notation}
$\Gamma_{ii'}$ stands for $\Gamma_i$ or $\Gamma_{i'}$. The relation
${}\langle\G_{a} , \G_{b}\rangle = e$ means $\G_{a} \G_{b} \G_{a} = \G_{b} \G_{a} \G_{b}$.
\end{notation}

\begin{theorem}\label{pres1}
The group  $\pi_1 (\C^2 \setminus B_1)$ is generated by
$\set{\G_j}_{j=1}^4$ subject to the relations
\begin{eqnarray}
{}\langle\G_{i} , \G_{i+1}\rangle & = & e \ \ \emph{for i=1,2,3} \label{sof1}\\
{}[ \G_{1} , \G_{i} ] & = & e \ \ \emph{for i=3,4}\label{sof6} \\
{}[ \G_{2} , \G_{4} ] & = & e \label{sof8}\\
{}[ \G_{4} , \G_{3} \G_{2} \G_{1}^2 \G_{2} \G_{3}] & = & e. \label{sof9}
\end{eqnarray}

The group $\pi_1 (\C\P^2 \setminus B_1)$ is isomorphic to
${\mathcal{B}_5}/{\langle\G_{4}^2  \G_{3} \G_{2} \G_{1}^2 \G_{2} \G_{3}\rangle}$,
and the group $\Pi_{(B_1)}$ is isomorphic to $S_5$.
\end{theorem}

\begin{proof}
The group $\pi_1 (\C^2 \setminus B_1)$ is generated by the elements
$\set{\G_j,\G_{j'}}_{j=1}^4$,  where $\G_j$ and
$\G_{j'}$ are loops in $\C^2$ around $j$ and $j'$, respectively.

By the \vK \ Theorem, the two branch points braids give the
following relations
\begin{eqnarray}
\G_i & = & \G_{i'} \ \ \mbox{for i=1,4}. \label{1}
\end{eqnarray}
From the monodromies $\varphi_{m_2}, \varphi_{m_4}$ and $\varphi_{m_5}$,
we produce relations (\ref{3})-(\ref{5}), (\ref{x1a})-(\ref{8}) and (\ref{x1b})-(\ref{11}) respectively
(e.g., from Figure \ref{3phi} we have (\ref{3})-(\ref{5})):
\begin{eqnarray}
\langle\G_{2'} , \G_{33'}\rangle = \langle\G_{2'} , \G_3^{-1}\G_{3'}\G_3\rangle
& = & e \label{3}\\
\G_{3'} \G_{3} \G_{2'} \G_{3}^{-1} \G_{3'}^{-1} & = & \G_2 \label{5}\\
\langle\G_{11'} , \G_{2}\rangle = \langle\G_1^{-1}\G_{1'}\G_1 , \G_2\rangle & = & e
\label{x1a}\\
\G_{2} \G_{1'} \G_{1} \G_{2}  \G_{1}^{-1} \G_{1'}^{-1} \G_{2}^{-1}
& = & \G_{2'} \label{8}\\
\langle\G_{33'} , \G_4\rangle = \langle\G_3^{-1}\G_{3'}\G_3 , \G_4\rangle & = & e
\label{x1b}\\
\G_{4} \G_{3'} \G_{3} \G_{4} \G_{3}^{-1} \G_{3'}^{-1} \G_{4}^{-1}
& = & \G_{4'}. \label{11}
\end{eqnarray}
The parasitic intersections braids contribute commutative relations
\begin{eqnarray}
{}[ \G_{11'} , \G_{ii'} ] & = & e \ \ \mbox{for i=3,4}\label{12} \\
{}[ \G_{22'} , \G_{44'} ] & = & e \label{14}.
\end{eqnarray}

Using (\ref{1}), (\ref{x1a}) and (\ref{x1b}), relations (\ref{8}) and (\ref{11})
can be rewritten as $\G_1^{-2} \G_2 \G_1^2 = \G_{2'}$ and $\G_4^{-2} \G_3 \G_4^2 = \G_{3'}$, respectively.
Using the fact that $\langle\G_{2} , \G_{3'}\rangle = \langle \G_{1}^2 \G_{2'} \G_{1}^{-2}, \G_{3'}\rangle =1$,
we can rewrite (\ref{5}) as  $\G_{2}^{-1} \G_{3} \G_{2'} \G_{3}^{-1} \G_{2}  =  \G_{3'}$.
Substituting these three relations in one another gives (\ref{sof9}), and substituting them in
(\ref{3}), (\ref{x1a}) and (\ref{x1b}) (in (\ref{12}) and (\ref{14}), respectively),
gives (\ref{sof1}) ((\ref{sof6}) and (\ref{sof8}), respectively).

In order to get $\pi_1 (\C\P^2 \setminus B_1)$, we add the projective relation
$\G_{4'} \G_{4} \G_{3'} \G_{3} \G_{2'} \G_{2} \G_{1'} \G_{1} =~e$, which is
transformed to $\G_{4}^2 \G_{3} \G_{2} \G_{1}^2 \G_{2} \G_{3} =e$. Therefore,
relation (\ref{sof9}) is omitted and  we have $\pi_1 (\C\P^2 \setminus B_1) \ \isom  \
{\mathcal{B}_5}/{\langle\G_{4}^2 \G_{3} \G_{2} \G_{1}^2 \G_{2} \G_{3}\rangle}$ and $\Pi_{(B_1)} \ \isom \ S_5$.
\end{proof}

\subsection{The surface $X_2$}
In \cite{MT1}, Moishezon-Teicher embed the surface $X_2=\C\P^1 \times \C\P^1$
into a big projective space by the linear system $(\mathcal{O}(i), \mathcal{O}(j))$,
where $i \geq 2, j \geq 3$. They use its degeneration to compute the fundamental
group of the Galois cover corresponding to the generic projection of the surface
onto $\C\P^2$.

\begin{figure}[h]
\begin{minipage}{\textwidth}
\begin{center}
\epsfbox{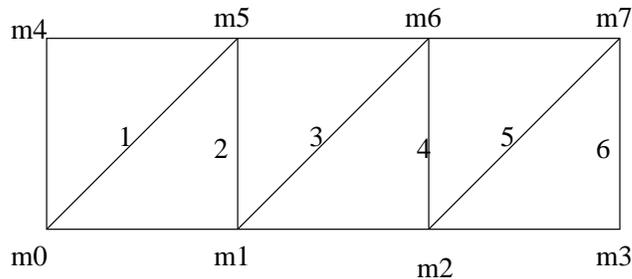} \caption{Degeneration of $X_2$}\label{1on3}
\end{center}
\end{minipage}
\end{figure}

In this paper the embedding is by the linear system $(\mathcal{O}(3), \mathcal{O}(1))$.
The degeneration of $X_2$ is a union of six planes embedded in $\C\P^7$, as depicted in Figure \ref{1on3}.

Now we explain what happens in the regeneration of the branch curve ${(B_2)}_0$.
Each diagonal line regenerates to a conic.
That means that in neighborhoods of $m_0$ and $m_7$ we have only conics, while in
neighborhoods of $m_1, m_2, m_5, m_6$  the conics are tangent to the lines they intersected with
(the vertical lines in the figure). Then each of these lines regenerates, causing a
regeneration of each tangency to three cusps. We end up with the curve $B_2$ with degree $10$
and $12$ cusps.

\begin{theorem}\label{del10}
The braid monodromy factorization of the curve $B_2$ is the product of
\begin{eqnarray}
\varphi_{m_0} & = & Z_{1 \; 1'}\\
\varphi_{m_1} & = & Z^3_{2 \; 2', 3} \cdot {Z_{3\; 3' }}^{Z^2_{2 \; 2', 3}}\\
\varphi_{m_2} & = & Z^3_{4 \; 4', 5} \cdot {Z_{5\; 5' }}^{Z^2_{4 \; 4', 5}}\\
\varphi_{m_5} & = & Z^3_{1', 2 \; 2'} \cdot {Z_{1\; 1'}}^{Z^2_{1', 2 \; 2'}}\\
\varphi_{m_6} & = & Z^3_{3', 4 \; 4'} \cdot {Z_{3\; 3'}}^{Z^2_{3', 4 \; 4'}}\\
\varphi_{m_7} & = & Z_{5 \; 5'}
\end{eqnarray}
and the parasitic intersections braids
\begin{eqnarray}\label{nelli}
{Z^2}_{1 \; 1',3 \; 3'},  \; {\stackrel{\scriptstyle (3)(3')}{Z^2}_{1 \; 1',4 \; 4'}},  \;
{\stackrel{\scriptstyle (3)(3')}{Z^2}_{2 \; 2',4 \; 4'}}, \; {Z^2}_{1 \; 1',5 \; 5'}, \;
{Z^2}_{2 \; 2',5 \; 5'}, \; {Z^2}_{3 \; 3',5 \; 5'}.
\end{eqnarray}
\end{theorem}

\begin{proof}
The monodromies $\varphi_{m_0}$ and $\varphi_{m_7}$ are braids of
branch points of the conics there. The monodromies $\varphi_{m_1}, \varphi_{m_2}$ (resp. $\varphi_{m_5},
\varphi_{m_6}$) are similar to the monodromies (\ref{m4}) and (\ref{m5}) (resp. (\ref{m2})).
According to this similarity of braids (only modification of indices in Figure \ref{3phi}),
we depict only the parasitic intersections braids in Figure \ref{ci1}.
\begin{figure}[h]
\begin{minipage}{\textwidth}
\begin{center}
\epsfbox{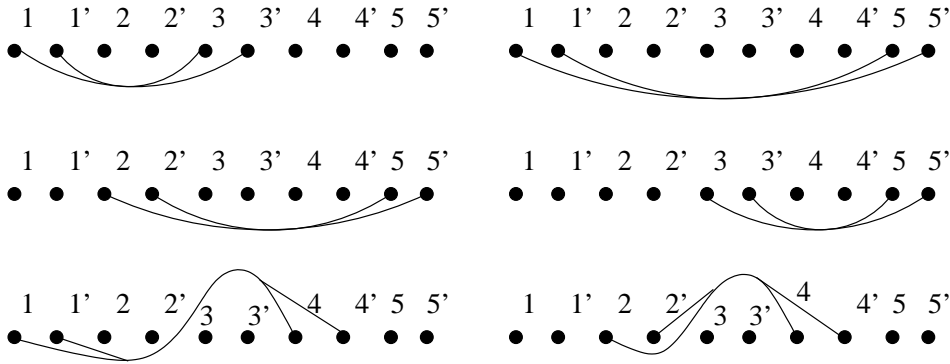} \caption{Parasitic intersections braids in the
factorization of $B_2$}\label{ci1}
\end{center}
\end{minipage}
\end{figure}
\end{proof}

We apply the van Kampen Theorem on the above braids to get a
presentation for $\pi_1(\C^2 \setminus B_2)$, and by omitting the generators
$\G_i, i= 1, \dots, 5$, and simplifying the relations, as done in the proof of
Theorem \ref{pres1}, we get the following:

\begin{theorem}\label{cp1on3}
The fundamental group $\pi_1(\C^2 \setminus B_2)$ is generated by
$\{\G_j\}_{j=1}^5$ subject to the relations
\begin{eqnarray}
{}\langle\G_{i} , \G_{i+1}\rangle & = & e  \ \ \emph{for i=1,2,3,4} \label{b3a}\\
{}[\G_{1} , \G_{i}] & = & e \ \ \emph{for i=3,4,5} \\
{}[\G_{2} , \G_{i}] & = & e \ \ \emph{for i=4,5} \\
{}[\G_{3} , \G_{5}] & = & e \\
{}\G_{2}^{-1} \G_{1}^{-2} \G_2^{-1} \G_{3} \G_2 \G_{1}^{2}
\G_2 & = & \G_{4}^{-1} \G_{5}^{-2} \G_{4}^{-1} \G_3 \G_{4} \G_{5}^2 \G_{4}. \label{b3b}
\end{eqnarray}

The group $\pi_1(\C\P^2 \setminus B_2)$ is isomorphic to ${\mathcal{B}_6}/
{\langle\G_{3} \G_4 \G_{5}^{2} \G_{4} \G_3 \G_{2} \G_1^2 \G_{2}\rangle}$, and the group
$\Pi_{(B_2)}$ is isomorphic to $S_6$.
\end{theorem}

\forget
\begin{proof}
By the van Kampen Theorem \cite{vK}, relation (\ref{new1}) comes from $\varphi_{m_0}$ and
$\varphi_{m_7}$, relations  (\ref{sofx4a}) -  (\ref{new14}) come from $\varphi_{m_1},
\varphi_{m_2}, \varphi_{m_5}, \varphi_{m_6}$, and relations (\ref{new15}) - (\ref{sofx4b})
come from (\ref{nelli}).
\begin{eqnarray}
{}\G_i & = & \G_{i'} \ \ \emph{for i=1,5} \label{new1}\\
{}\langle\G_{22'},\G_3\rangle = \langle\G_2^{-1}\G_{2'}\G_2,\G_3\rangle & = & e \label{sofx4a}\\
{}\G_{3} \G_{2'} \G_2 \G_{3} \G_2^{-1} \G_{2'}^{-1}
\G_3^{-1} & = & \G_{3'} \label{new5}\\
{}\langle\G_{44'} , \G_{5}\rangle = \langle\G_4^{-1} \G_{4'} \G_4 , \G_5\rangle & = & e \\
{}\G_{5} \G_{4'} \G_{4} \G_{5} \G_{4}^{-1} \G_{4'}^{-1} \G_{5}^{-1}
& = & \G_{5'} \label{new8}\\
{}\langle\G_{1'},\G_{22'}\rangle = \langle\G_{1'},\G_2^{-1}\G_{2'}\G_2 \rangle & = & e \\
{}\G_{2'} \G_2 \G_{1'} \G_{2}^{-1} \G_{2'}^{-1} & = & \G_{1} \label{new11}\\
{}\langle\G_{3'},\G_{44'}\rangle = \langle\G_{3'},\G_4^{-1}\G_{4'}\G_4 \rangle & = & e \\
{}\G_{4'} \G_4 \G_{3'} \G_{4}^{-1} \G_{4'}^{-1} & = & \G_{3} \label{new14}\\
{}[ \G_{11'} , \G_{ii'} ] & = & e \ \ \emph{for i=3,5} \label{new15} \\
{}[ \G_{ii'} , \G_{55'} ] & = & e \ \ \emph{for i=2,3} \label{new17} \\
{}[ \G_{3'} \G_{3} \G_{ii'} \G_{3}^{-1} \G_{3'}^{-1} , \G_{4} ] & = & e \ \ \emph{for i=1,2}\\
{}[ \G_{3'} \G_{3} \G_{ii'} \G_{3}^{-1} \G_{3'}^{-1} , \G_{4}^{-1}\G_{4'}\G_4 ] & = & e \ \
\emph{for i=1,2}, \label{sofx4b}
\end{eqnarray}
and a long process of simplification gives (\ref{b3a}) - (\ref{b3b}).

The projective relation $\G_{5'} \G_5 \cdots \G_{1'} \G_1 = e$ is transformed to
$\G_{3} \G_{4} \G_{5}^2 \G_{4} \G_{3} \G_2 \G_1^2 \G_2 =e$. This one cancels
(\ref{b3b}), and therefore $\pi_1(\C\P^2 \setminus B_2) \isom
{\mathcal{B}_6}/{<\G_{3} \G_4 \G_{5}^{2} \G_{4} \G_3 \G_{2} \G_1^2 \G_{2}>}$
and $\Pi_{(B_2)} \isom S_6$.
\end{proof}
\forgotten

One can easily generalize this result. Take $X_2:=\C\P^1\times \C\P^1$ embedded in
$\C\P^{n+1}$ by ${\mathcal O}(n,1)$. This embedded toric surface corresponds to the
convex polygon \linebreak $P_{n,1}:=\mbox{Conv}\{(0,0), (n,0), (0,1), (n,1)\}$.
And we have

\begin{corollary}
The groups $\Pi_{(B)}$ and $C$ are isomorphic to $S_{2n}$, and $\pi_1((X_2)_{Gal})$
is trivial.
\end{corollary}

\subsection{The surface $X_3$}
\begin{figure}[h]
\begin{minipage}{\textwidth}
\begin{center}
\epsfbox{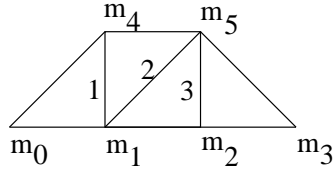}
\end{center}
\end{minipage}
\caption{Degeneration of $X_3$}\label{trapez}
\end{figure}

The degeneration of $X_3$ is a union of four planes embedded in $\C\P^5$,
as depicted in Figure \ref{trapez}.

The branch curve ${(B_3)}_0$ in $\C\P^2$ is a line arrangement. Regenerating it,
the diagonal line regenerates to a conic, which is tangent to the lines $1$ and $3$.
When the lines regenerate, each tangency regenerates into three cusps.
We obtain the branch curve $B_3$, whose degree is $6$ and which has six cusps.

\begin{theorem}
The braid monodromy factorization related to $B_3$ is the product of
\begin{eqnarray}
\varphi_{m_1} & = & Z^3_{1 \; 1', 2} \cdot {Z_{2\; 2' }}^{Z^2_{1\;1', 2}} \\
\varphi_{m_5} & = & Z^3_{2', 3\;3'} \cdot {Z_{2\; 2'}}^{Z^2_{2', 3\; 3'}} \\
\varphi_{m_2} & = & Z_{3\; 3'} \\
\varphi_{m_4} & = & Z_{1 \; 1'}
\end{eqnarray}
and the parasitic intersections braids
\begin{eqnarray}
{Z^2}_{\hspace{-.1cm}{1 \; 1',3 \; 3'}}.
\end{eqnarray}
\end{theorem}

\begin{proof}
Similar proof as in Theorem \ref{del8}.
\end{proof}

We apply the van Kampen Theorem on the above braids to get a
presentation for $\pi_1(\C^2 \setminus B_3)$, and again, by simplifying
the relations and omitting generators, we get:

\begin{theorem}\label{pretrapez}
The fundamental group  $\pi_1(\C^2 \setminus B_3)$ is generated by
$\G_1, \G_{2}, \G_3$ subject to the relations
\begin{eqnarray}
{}\langle\G_{i} , \G_{i+1}\rangle & = & e \ \ \emph{for i=1,2} \label{x2a}\\
{}[ \G_{1} , \G_{3} ] & = & e\label{b2co}\\
{}\G_1^{-2} \G_{2} \G_1^{2} & = & \G_3^{-2} \G_{2} \G_3^{2}\label{tom}.
\end{eqnarray}

The group $\pi_1(\C\P^2 \setminus B_3)$ is isomorphic to ${\mathcal{B}_4} / {\langle\G_2
\G_3^2 \G_2 \G_1^2\rangle}$, and the group $\Pi_{(B_3)}$ is isomorphic to $S_4$.
\end{theorem}

\forget
\begin{proof}
By the van Kampen Theorem \cite{vK},
\begin{eqnarray}
\ \G_i & = & \G_{i'} \ \ \emph{for i=1,3} \label{newer1}\\
\ \langle\G_{11'} , \G_{2}\rangle = \langle\G_1^{-1}\G_{1'}\G_1, \G_2\rangle & = & e \\
\ \G_{2} \G_{1'} \G_1 \G_{2} \G_1^{-1} \G_{1'}^{-1}
\G_2^{-1} & = & \G_{2'} \label{newer5}\\
\ \langle\G_{2'} , \G_{33'}\rangle = \langle\G_{2'} , \G_3^{-1} \G_{3'} \G_3\rangle & = & e \\
\ \G_{3'} \G_{3} \G_{2'} \G_{3}^{-1} \G_{3'}^{-1}
& = & \G_{2} \label{newer8}\\
\ [ \G_{11'} , \G_{33'} ] & = & e. \label{newer9}
\end{eqnarray}
From this presentation one can easily derive relations (\ref{x2a}) - (\ref{tom}).

The projective relation $\G_{3'} \G_3 \G_{2'} \G_2 \G_{1'} \G_1 = e$ gets the form
$\G_2 \G_3^2 \G_2 \G_1^2=e$. Substituting it in (\ref{tom}), we get again
$\G_2 \G_3^2 \G_2 \G_1^2=e$. Therefore we are left with (\ref{x2a}), (\ref{b2co})
and with the above relation. Hence $\pi_1(\C\P^2 \setminus B_3)
\isom \mathcal{B}_4/ <\G_2 \G_3^2 \G_2 \G_1^2>$ and $\Pi_{(B_3)} \isom S_4$.
\end{proof}
\forgotten

\section{The surface $X_4$}\label{x4}
\begin{figure}[h]
\begin{minipage}{\textwidth}
\begin{center}
\epsfbox{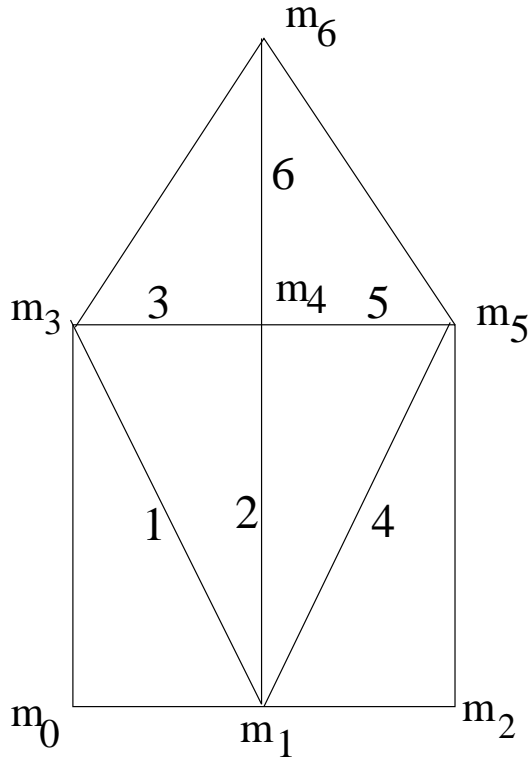}
\end{center}
\end{minipage}
\caption{Degeneration of $X_4$}\label{home}
\end{figure}

The degeneration ${(X_4)}_0$ of $X_4$ is a union of six planes embedded in $\C\P^6$
(Figure \ref{home}).

\begin{figure}[h]
\begin{minipage}{\textwidth}
\begin{center}
\epsfbox{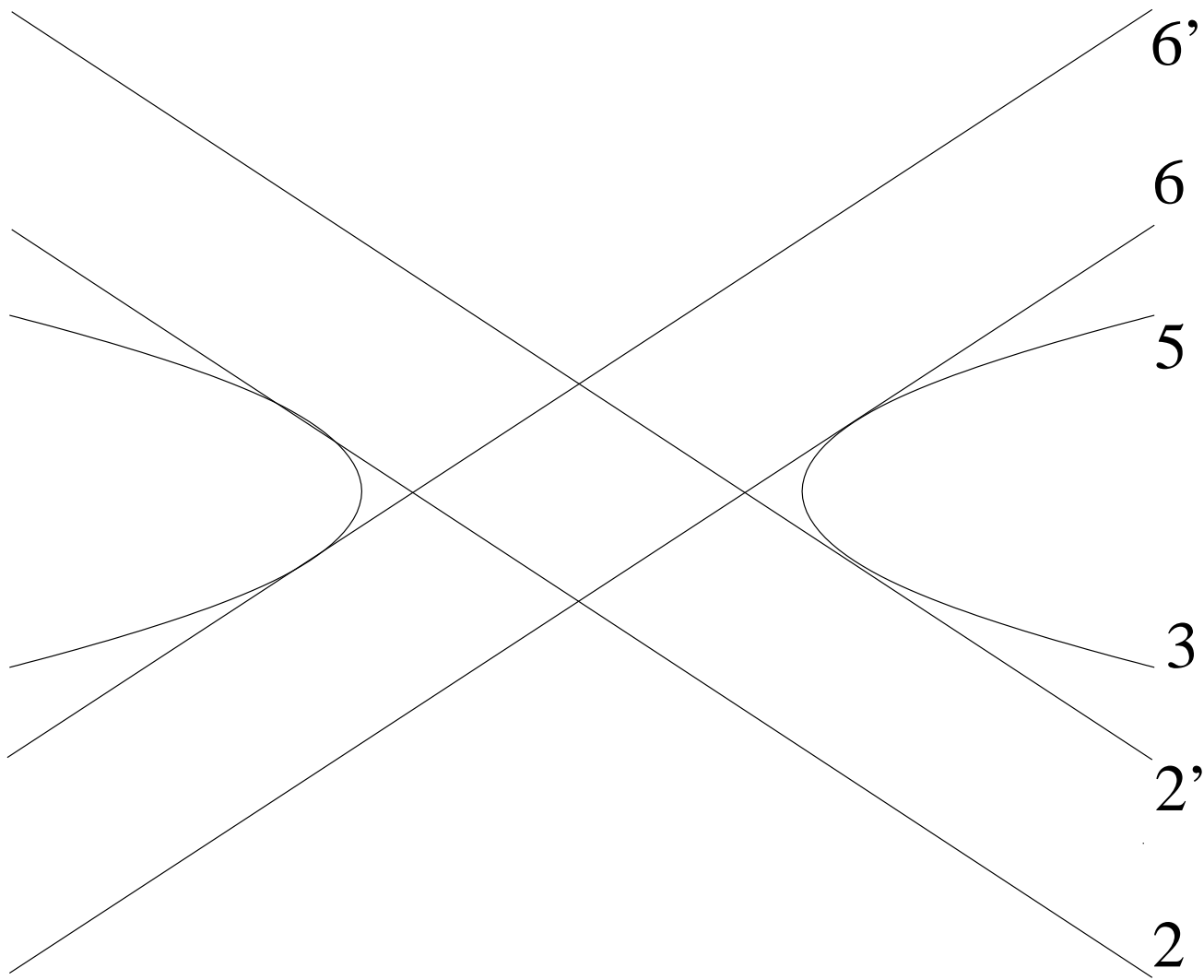}
\end{center}
\end{minipage}
\caption{Regeneration around the $4$-point $m_4$}\label{4lines}
\end{figure}

The regeneration of ${(X_4)}_0$ induces a regeneration on the branch curve
${(B_4)}_0$ (line arrangement, composed of six lines).
$X_4$ has $A_1$ singularity as explained in the introduction.
That means that the regeneration of the top vertex $m_6$ should yield a node in
the branch curve, involving the components labelled 6 and 6' (so that
the double cover possesses an ordinary double point). The vertices
$m_3$ and $m_5$ are $2$-points, and therefore the regeneration around
them is already known: the line $1$ ($4$ resp.) regenerates to a conic
which is tangent to the line $3$ ($5$ resp.). When these lines regenerate,
each tangency regenerates to three cusps. The vertex $m_4$ is a $4$-point,
see e.g., \cite{AT-TT1}. The regeneration is as follows.
The lines $3$ and $5$ regenerate to a hyperbola, and each line among $2$ and $6$
regenerates to a pair of parallel lines. The hyperbola is then tangent to the
lines $2, 2', 6, 6'$, see Figure \ref{4lines}.
The hyperbola doubles, therefore we have four branch points, and moreover,
each tangency regenerates to three cusps.

\begin{figure}[h]
\begin{minipage}{\textwidth}
\begin{center}
\epsfbox{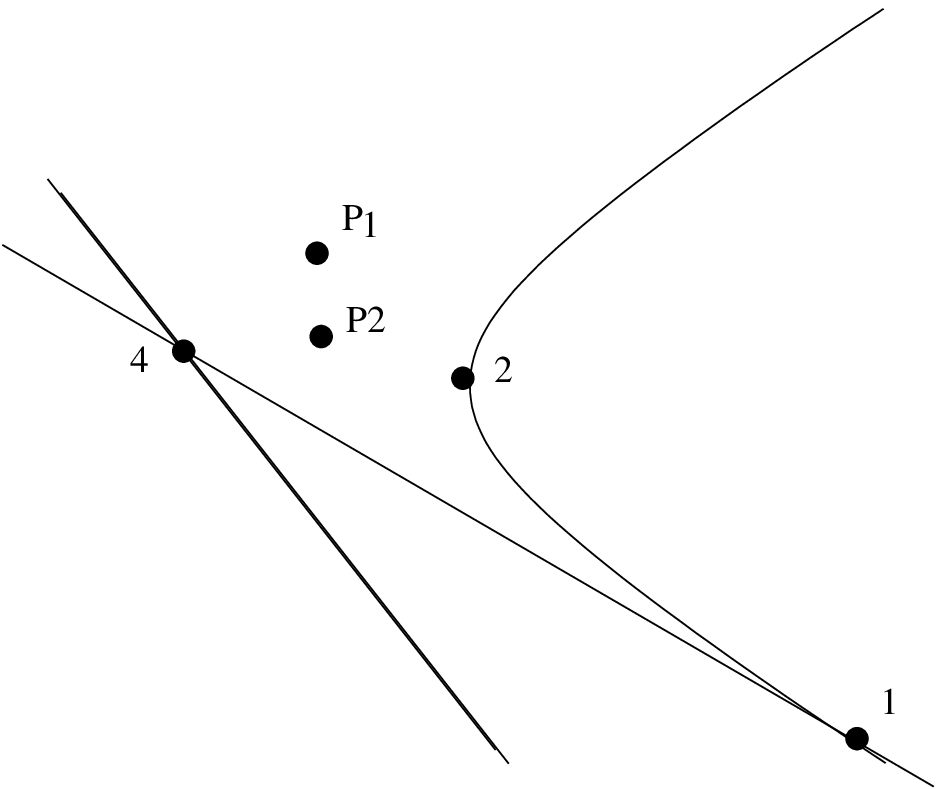}
\end{center}
\end{minipage}
\caption{Regeneration around $m_1$}\label{fignew}
\end{figure}

However, the vertex $m_1$
is of new type. The regeneration can be done as follows. Line $4$ regenerates to a conic,
while $1$ is still unregenerated. Figure \ref{fignew} describes this step. The points $P_1$ and $P_2$
are the intersections of $1$ with the conic (they are complex).
The intersection of lines $1$ and $2$
can be then locally considered as a $2$-point; this means that $1$ regenerates to a conic, which is
tangent to line $2$. At this point $P_1$ and $P_2$ are doubled.
Line $2$ then regenerates to a pair of parallel lines $2$ and $2'$,
and each tangency regenerates to three cusps.
Note that keeping a parabola, which we get in the regeneration
around $m_1$, as our picture in the affine part of the conics,
we have possibly another branch point further away, possibly at infinity.
We prove below the existence of these two extra branch points,
which contribute two half-twists to the braid monodromy factorization.

The parasitic intersections are fixed by Figure \ref{home} and this time they are
the intersections in $\C\P^2$ of line $1$ with lines $5$ and $6$, and $4$ with $3$ and $6$.

Therefore, we have

\begin{theorem}\label{bmtrapez}
The braid monodromies which we get from the regeneration around ${m_1}, m_3, m_4, m_5, m_6$
are
\begin{eqnarray}
{}\varphi_{m_1} & = & Z^3_{2 \; 2', 4} \cdot {(Z_{4 \; 4'})}^{Z^2_{2 \; 2', 4}}
\cdot (Z^2_{1 \; 1', 4'} \cdot {(Z^2_{1 \; 1', 4})}^{Z^2_{2 \; 2', 4}}) \cdot
({Z^3_{1', 2 \; 2'}  \cdot (Z_{1 \; 1'})}^{Z_{1', 2 \; 2'}^2}) \\
{}\varphi_{m_3} & = & Z^3_{1', 3 \; 3'} \cdot {(Z_{1 \; 1'})}^{Z^2_{1', 3 \; 3'}}\\
{}\varphi_{m_5} & = & Z^3_{4', 5 \; 5'} \cdot {(Z_{4 \; 4'})}^{Z^2_{4', 5 \; 5'}}\\
{}\varphi_{m_6} & = & Z^2_{6 \; 6'} \\
{}\varphi_{m_4} & = & (Z^3_{2', 3 \; 3'} \cdot {Z^3_{5 \; 5',6}} \cdot h_1 \cdot h_2
\cdot {(Z^2_{2' \; 6})}^{Z_{2', 3 \; 3'}^2} \cdot Z^2_{2 \; 6}) \cdot \\ & &
(Z^3_{2, 3 \; 3'} \cdot {(Z^3_{5 \; 5',6'})}^{Z^{-2}_{6 \; 6'}}
\cdot h_3 \cdot h_4 \cdot {(Z^2_{2 \; 6'})}^{Z_{2, 3 \; 3'}^2 Z^{-2}_{6 \; 6'}} \cdot
{(Z^2_{2' \; 6'})}^{Z_{2 \; 2'}^{-2} Z^{-2}_{6 \; 6'}} ),\nonumber
\end{eqnarray}
where  $h_1, h_2$ are the upper braids and $h_3, h_4$ are the lower ones in Figure \ref{h1234}.
\begin{figure}[h]
\begin{minipage}{\textwidth}
\begin{center}
\epsfbox{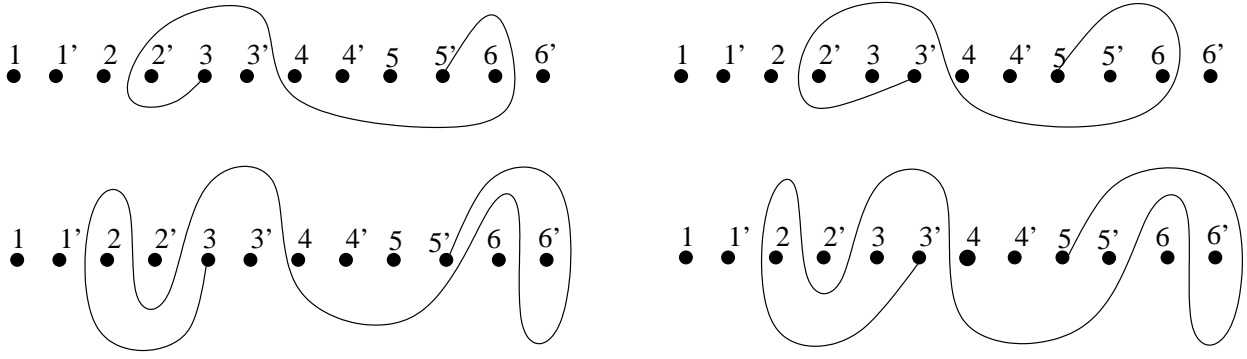}
\end{center}
\end{minipage}
\caption{The braids $h_1, h_2, h_3, h_4$}\label{h1234}
\end{figure}

The parasitic intersections braids (Figure \ref{parax4}) are
\begin{equation}
{(Z^2_{1 \; 1', 5 \; 5'})}^{Z^{-2}_{4 \; 4', 5 \; 5'}}, \;
Z^2_{1 \; 1', 6 \; 6'}, \;  Z^2_{3 \; 3', 4 \; 4'}, \; Z^2_{4 \; 4', 6 \; 6'}.
\end{equation}
\begin{figure}[h]
\begin{minipage}{\textwidth}
\begin{center}
\epsfbox{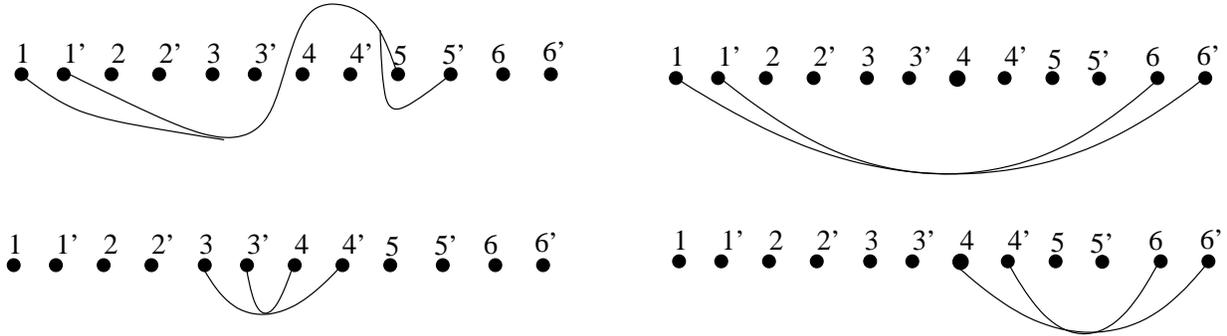}
\end{center}
\end{minipage}
\caption{Parasitic intersections braids in the factorization of $B_4$}\label{parax4}
\end{figure}
\end{theorem}

Since $B_4$ has  degree $12$, the total degree of the braid monodromy factorization $\Delta_{12}^2$
should be $12 \cdot 11 =132$, see \cite{MT2}.
By the above regeneration, $B_4$ has $8$ branch points, $24$ cusps, and $25$ nodes.
Their related braids give a total degree of $130$. The missing braids correspond to
two extra branch points. We explain how to find them.

We look at the preimage in $X_4$ of a vertical line in $\C\P^2$
(a fiber of the projection); this is an elliptic curve (a $6$-fold cover of $\C\P^1$ branched in
$12$ points). Considering the entire family of vertical lines in $\C\P^2$, we get
that $X_4$ admits a projection to $\C\P^1$, with generic fiber an elliptic curve.
The preimage of a vertical line in $\C\P^2$ is singular if and only if that
vertical line is tangent to the branch curve or if it passes through the intersection of
the lines $6$ and $6'$.

There is a ``lifting homomorphism" from the braid group
$B_{12}$ to the mapping class group $SL(2,\Z)$, obtained by considering the above-mentioned
$6$-fold cover of $\C\P^1$: if the 12 branch points are moved by a braid,
this induces a homeomorphism of the covering, \cite[Section 5.2]{AK}.
Now, since the abelianization of $SL(2,Z)$ is $\Z/12$ and the quotient homomorphism
$SL(2,Z) \rightarrow \Z/12$ takes Dehn twists to the integer 1,
the number of Dehn twists we get is a multiple of 12.
However, we get $2$ from $Z^2_{6 \; 6'}$, and $1$ from each one of the $8$ branch points.

In order to check which braids are missing, we consider a homomorphism from the pure
braid group on $12$ strings to the pure braid group on $2$ strings, defined
by deleting all the strands except $i$ and $i'$; it should map $\Delta_{12}^2$
to $\Delta_{2}^2 = Z_{i \; i'}^2$. By Lemma~2.I in \cite{MT5}, $Z^3_{i \; i', j} =
Z^2_{i' \; j} Z^2_{i \; j} Z^2_{i' \; j} Z^2_{i \; j} Z_{i \; i'}$.
Therefore, by Theorem \ref{bmtrapez}, we get $\Delta_{2}^2 = Z_{i\; i'}^2$ for
$i = 1, 2, 4, 6$. Now, forgetting all indices and remembering $3$ and $3'$ (resp. $5$ and $5'$)
gives the half-twist $Z_{3 \; 3'}$  (resp. $Z_{5 \; 5'}$), counted three times.
But by Lemma 8.IV in \cite{MT5}, $\varphi_{m_4} =
\Delta_8^{2} Z^{-2}_{2 \; 2'} Z^{-2}_{6 \; 6'} Z^{-2}_{3 \; 3'} Z^{-2}_{5 \; 5'}$.
In his thesis \cite{robb}, Robb discusses existence of extra branch points. According to
our results, there is an extra branch point, which contributes the half-twist $Z_{3 \; 3'}$
(resp. $Z_{5 \; 5'}$). By \cite[Prop. 3.3.1]{robb}, the relation in
$\pi_1(\C\P^2 \setminus B_4)$ should be $\G_{3} = \G_{3'}$
(resp. $\G_{5} = \G_{5'}$).
\begin{remark}
Another  justification for this can be also group-theoretic.
Because Moishezon-Teicher's
formulas for arrangements of lines \cite{MT5} deal only with what happens before each
line regenerates to a pair $i, i'$, their global formula $($$\Delta^2 = \prodl C_i \varphi_i$,
$C_i$ are the parasitic braids$)$
is only correct up to half-twists of the form $Z_{i \; i'}$, which are not seen at
all by configurations at the level of the double lines $($before regeneration$)$.
In our case the above product is not $\Delta_{12}^2$ but
$\Delta_{12}^2 Z_{3 \; 3'}^{-1} Z_{5 \; 5'}^{-1}$, and thus implies that there are two
extra half-twists which must be $Z_{3 \; 3'}$ and $Z_{5 \; 5'}$.
\end{remark}

\begin{corollary}
The braid monodromy factorization $\Delta_{12}^2$ is a product of the braids from
Theorem \ref{bmtrapez} and the extra branch points braids $Z_{3 \; 3'}$ and $Z_{5 \; 5'}$.
\end{corollary}

Now we are ready to compute the group $\pi_1(\C\P^2 \setminus B_4)$.

\begin{theorem}\label{prehome}
Let $\tilde{\B}_6$ be the quotient of the braid group $\B_6$ by $\langle[X,Y]\rangle$, where $X, Y$ are transversal.
The fundamental group $\pi_1(\C\P^2 \setminus B_4)$ is isomorphic to a quotient of $\tilde{\B}_6$
by $\langle(\ref{anilo})\rangle$. The group  $\Pi_{(B_4)}$ is isomorphic to $S_6$.
\end{theorem}

\begin{proof}
Applying the van Kampen Theorem \cite{vK} on the factorization $\Delta_{12}^2$, we get a
presentation of $\pi_1(\C\P^2 \setminus B_4)$ with the generators $\set{\G_i, \G_{i'}}_{i=1}^6$.

The monodromy $\varphi_{m_1}$ contributes the relations
\begin{eqnarray}
{}\langle\G_{22'} , \G_{4}\rangle = \langle\G_{2'} \G_{2} \G_{2'}^{-1} , \G_4\rangle & = & e \label{tr23} \\
{}\G_{4}^{\G_{2}^{-1} \G_{2'}^{-1} \G_{4}^{-1}} & = & \G_{4'}\label{br23}\\
{}[ \G_{11'}, \G_{4}^{\G_{2}^{-1} \G_{2'}^{-1} \G_{4}^{-1}} ] = [\G_{11'}, \G_{4'}] & = & e \label{com14}\\
{}\langle\G_{1'} , \G_{22'}\rangle = \langle\G_{1'} , \G_{2'} \G_{2} \G_{2'}^{-1}\rangle & = & e \label{tr12} \\
{}\G_{2'}\G_{2}\G_{1'} \G_{2}^{-1} \G_{2'}^{-1} & = & \G_{1} \label{br12}.
\end{eqnarray}
From the monodromies $\varphi_{m_3}$ and $\varphi_{m_5}$ we have
\begin{eqnarray}
{}\langle\G_{1'} , \G_{33'}\rangle = \langle\G_{1'} , \G_{3'} \G_{3} \G_{3'}^{-1}\rangle & = & e\label{tr14}\\
{}\G_{3'} \G_{3} \G_{1'} \G_{3}^{-1} \G_{3'}^{-1} & = & \G_{1} \label{br14}\\
{}\langle\G_{4'} , \G_{55'}\rangle = \langle\G_{4'} , \G_{5'} \G_{5} \G_{5'}^{-1}\rangle & = & e \label{tr35}\\
{}\G_{5'} \G_{5} \G_{4'} \G_5^{-1} \G_{5'}^{-1} & = & \G_{4}. \label{br35}
\end{eqnarray}
By $\varphi_{m_4}$ we have
\begin{eqnarray}
{}\langle\G_{22'}, \G_{3}\rangle = \langle\G_{22'}, \G_{3'}\rangle =
  \langle\G_{22'}, \G_{3'} \G_{3} \G_{3'}^{-1}\rangle & = & e \label{tr24}\\
{}\langle\G_{55'}, \G_{6}\rangle = \langle\G_{5'} \G_{5} \G_{5'}^{-1} ,\G_{6}\rangle & = & e \label{tr56}\\
{}\langle\G_{55'}, \G_{6}^{-1} \G_{6'} \G_{6}\rangle = \langle\G_{5'} \G_{5} \G_{5'}^{-1},
  \G_{6}^{-1} \G_{6'} \G_{6}\rangle & = & e \label{tr56a}\\
{}[ \G_{2}, \G_{6} ] = [ \G_{2'}^{\G_2}, \G_{6'}^{\G_6} ] & = & e \label{com26a}\\
{}[ \G_{2'}, \G_{6}^{\G_{3'} \G_{3}} ] = [ \G_{2}, \G_{6'}^{\G_6 \G_{3'} \G_{3}} ]& = & e \label{com26b}\\
{}\G_{3}^{\G_{2'}^{-1} \G_{3}^{-1} \G_{3'}^{-1}} & = & \G_{5'}^{\G_6^{-1}} \label{4pta}\\
{}\G_{3'}^{\G_{2'}^{-1} \G_{3}^{-1} \G_{3'}^{-1}} & = & \G_{5}^{\G_{5'}^{-1} \G_6^{-1}} \label{4pta1}\\
{}\G_{3}^{\G_{2}^{-1} \G_{3}^{-1} \G_{3'}^{-1}} & = & \G_{5'}^{\G_{6}^{-1} \G_{6'}^{-1} \G_6} \label{4pta2}\\
{}\G_{3'}^{\G_{2}^{-1} \G_{3}^{-1} \G_{3'}^{-1}} & = & \G_{5}^{\G_{5'}^{-1} \G_{6}^{-1} \G_{6'}^{-1} \G_6},\label{4ptb}
\end{eqnarray}
and $\varphi_{m_6}$ contributes
\begin{eqnarray}
{}[ \G_{6}, \G_{6'} ] & = & e. \label{66'}
\end{eqnarray}
From the parasitic intersections braids, we have
\begin{eqnarray}
{}[ \G_{11'}, \G_{55'}^{\G_{4'} \G_{4}}] & = & e \label{com15}\\
{}[ \G_{11'} , \G_{66'} ] & = & e \label{com16}\\
{}[ \G_{44'} , \G_{ii'}] & = & e \ \ \mbox{for i=3,6}, \label{com346}
\end{eqnarray}
and the extra branch points contribute
\begin{eqnarray}
{} \G_{3} &=& \G_{3'} \label{44'}\\
{} \G_{5} &=& \G_{5'}. \label{55'}
\end{eqnarray}
The projective relation is
\begin{eqnarray}
{}\G_{6'}\G_{6}\G_{5'}\G_{5}\G_{4'}\G_{4}\G_{3'}\G_{3}\G_{2'}\G_{2}\G_{1'}\G_1 & = & e. \label{proj4a}
\end{eqnarray}

\begin{lemma}
The above presentation is a complete one.
\end{lemma}
{\it{Proof of the Lemma.}} Considering the complex conjugations (details in \cite{MT1}, \cite{MT5})
of the braids, we get a complete set of relations. Simplifying them gives the same list as above.

\bigskip

We outline now the simplification of the above presentation. We will express the relations
in terms of $\G_{1}, \G_2, \G_3, \G_4, \G_5$ and $\G_{6'}$. First we use relations (\ref{44'})
and (\ref{55'}) to omit the generators $\G_{3'}$ and $\G_{5'}$ from all the given relations.

The branch points relations (\ref{br23}), (\ref{br12}), (\ref{br14}), (\ref{br35}) and
(\ref{4pta}) - (\ref{4ptb}) are rewritten as
\begin{eqnarray}
{}\G_{4'} &=& \G_{2}^{-1} \G_{2'}^{-1} \G_4 \G_{2'} \G_2 \ \ \mbox{(using (\ref{tr23}))}\label{t4}\\
{}\G_{1'} &=& \G_{2}^{-1} \G_{2'}^{-1}\G_{1} \G_{2'} \G_{2} \label{t5}\\
{}\G_{1'} &=& \G_3^{-2} \G_1 \G_3^2\label{t1}\\
{}\G_{4'} &=& \G_{5}^{-2} \G_4 \G_5^{2}\label{t2}\\
{}\G_{6} &=& \G_{5} \G_3 \G_{2'} \G_{3}^{-1} \G_5^{-1}\label{t3}\\
{}\G_{6'} &=& \G_{5} \G_3 \G_2 \G_{3}^{-1} \G_5^{-1}.\label{t3a}
\end{eqnarray}

Now we rewrite the commutations. Using (\ref{t3a}), relation
(\ref{tr56}) gets the form $\langle\G_{3}, \G_{6}\rangle =e$, and
this enables us to prove that (\ref{com26b}) is

\begin{eqnarray}
{}e=[ \G_{2'} , \G_{3}^{-2} \G_{6} \G_{3}^{2}]=  [ \G_{3}^{-1}
\G_{6} \G_{5} \G_{6}^{-1} \G_3,
\G_{3}^{-2} \G_{6} \G_{3}^{2}] =  & & \label{35com}\\
{}[\G_{6} \G_{5} \G_{6}^{-1}, \G_{3}^{-1} \G_{6} \G_{3}] = [\G_{6}
\G_{5} \G_{6}^{-1}, \G_6 \G_3 \G_{6}^{-1}] = [ \G_{3} , \G_{5}].&
& \nonumber
\end{eqnarray}

Relation (\ref{com15}) is rewritten as $[\G_{1}, \G_{5}]=e$, using
(\ref{tr35}), (\ref{t2}), (\ref{t1}) and $[ \G_{3} , \G_{5}]=e$.
This enables us to prove from (\ref{com14}) that $[\G_{1},
\G_{4}]=e$.

Using these two resulting relations, together with (\ref{t3}),
(\ref{t4}) and (\ref{com346}), relation $[ \G_{1} , \G_{6} ]=e$ is
rewritten as follows:
\begin{eqnarray}
{}& &e=[ \G_{1} , \G_{6} ]= [ \G_{1} , \G_{5} \G_{3} \G_{2'} \G_{3}^{-1} \G_{5}^{-1}]=
[ \G_{1} , \G_{3} \G_{2'} \G_{3}^{-1}]= [\G_{3}^{-1} \G_{1} \G_{3}, \G_{2'}] =  \nonumber \\
{}& &[\G_{3}^{-1} \G_{1} \G_{3}, \G_{4}^{-1}  \G_{2} \G_{4'} \G_{2}^{-1} \G_4]=
[\G_{3}^{-1} \G_{1} \G_{3}, \G_{4'}^{-1} \G_{2} \G_{4'}]=
[\G_{3}^{-1} \G_{1} \G_{3}, \G_{2}]=[\G_{1}, \G_{3} \G_{2} \G_{3}^{-1}]. \nonumber
\end{eqnarray}
In a similar way, $[ \G_{1'} , \G_{6} ] = e$ can be rewritten as
$[ \G_{1} , \G_{3}^{-1} \G_{2} \G_{3}]=e$. Using (\ref{t2}) and $[
\G_{3} , \G_{5}]=e$, relation $[ \G_{3} , \G_{4'} ] = e$ gets the
form $[ \G_{3} , \G_{4} ] = e$.  Relation (\ref{tr56}) is
rewritten as

\begin{eqnarray}
{}e=\langle\G_{5}, \G_{6}\rangle=\langle \G_{5} , \G_{5} \G_{3}
\G_{2'} \G_{3}^{-1} \G_{5}^{-1}\rangle =
\langle\G_{5}, \G_{2'}\rangle=& & \label{25tr} \\
{}\langle\G_5, \G_{1}^{-1} \G_{2} \G_{1'} \G_{2}^{-1}
\G_{1}\rangle = \langle\G_5, \G_{1'}^{-1} \G_{2} \G_{1'}\rangle =
\langle\G_{5}, \G_{2}\rangle.& & \nonumber
\end{eqnarray}
Thus $[ \G_{44'} , \G_{6} ]=e$  are rewritten as
\begin{eqnarray}
{}[ \G_{4} , \G_{5} \G_{2} \G_{5}^{-1}] = [ \G_{4} , \G_{5}^{-1} \G_{2} \G_{5}] =e. \label{425}
\end{eqnarray}

Now, relation $[ \G_{2'} , \G_{6'} ] = e$ gets the form $[ \G_{2}
, \G_{6'} ] = e$, using (\ref{com16}), (\ref{tr12}) and
(\ref{t5}).  This relation, together with (\ref{t4}) and
(\ref{t5}) enable us to prove  that $[ \G_{1'} , \G_{6'} ] = e$
and $[ \G_{4'} , \G_{6'} ]=e$ get the forms $[ \G_{1} , \G_{6'} ]
= e$ and $[ \G_{4} , \G_{6'} ]=e$ respectively. Since
(\ref{tr56a}) can be rewritten as  $\langle\G_{3}, \G_{6'}\rangle
=e$, relation $[\G_{2}, \G_{3}^{-2} \G_{6'} \G_{3}^2]=e$ from
(\ref{com26b}) gets the form $[ \G_{3} , \G_{5}]=e$.

Relation $[ \G_{2}, \G_{6} ]=e$ from (\ref{com26a}) gets the form
\begin{eqnarray*}
{}& &e=[\G_{2}, \G_{6}]\stackrel{(\ref{t3})}{=}[\G_{2} , \G_{5}
\G_{3} \G_{2'} \G_{3}^{-1} \G_{5}^{-1}]= [\G_{5}^{-1} \G_{2}
\G_{5}, \G_{3} \G_{2'} \G_{3}^{-1}]\stackrel{(\ref{t4})}{=}\\ & &
[\G_{5}^{-1} \G_{2} \G_{5}, \G_{3} \G_{4}^{-1}\G_{2} \G_{4'}
\G_2^{-1} \G_4 \G_{3}^{-1}] \stackrel{(\ref{com346}), (\ref{425}),
(\ref{25tr})}{=} [\G_{2}^{-1} \G_{3}^{-1} \G_{2} \G_{5}
\G_{2}^{-1}\G_{3} \G_{2},  \G_{4'}] \stackrel{(\ref{tr24})}{=}\\ &
& [\G_{3}\G_2^{-1} \G_{3}^{-1} \G_{5} \G_3\G_2\G_{3}^{-1},
\G_{4'}] \stackrel{(\ref{com346}),
(\ref{35com})}{=}[\G_2^{-1}\G_{5}\G_2, \G_{4'}]
\stackrel{(\ref{t4})}{=}\\& &[\G_2^{-1}\G_{5}\G_2,\G_{2}^{-1}
\G_{2'}^{-1} \G_4 \G_{2'} \G_2]
\stackrel{(\ref{tr23})}{=}[\G_{5},\G_4 \G_{2'}\G_4^{-1}]
\stackrel{(\ref{t5})}{=}\\& &[\G_{5},\G_4 \G_{1}^{-1} \G_2 \G_{1'}
\G_2^{-1}\G_{1}\G_4^{-1}] \stackrel{(\ref{com14}), (\ref{com15}),
(\ref{tr12})}{=} [\G_{5},\G_4 \G_{1'}^{-1} \G_2
\G_{1'}\G_4^{-1}]=\\& &[\G_{5},\G_4  \G_2 \G_4^{-1}]
\stackrel{(\ref{tr23}), (\ref{25tr})}{=}[\G_{4},\G_5^{-1}  \G_2
\G_5].
\end{eqnarray*}
The only relation which is left for now in its original form is
(\ref{66'}). We prove below that $\G_{6}= \G_{6'}$, and this
equality will eliminate it.

The triple relations are rewritten as follows. (\ref{tr14}) and
(\ref{tr35}) get the forms $\langle\G_{1}, \G_{3}\rangle =e$ and
$\langle\G_{4}, \G_{5}\rangle  = e$, using (\ref{t1}) and
(\ref{t2}). It is also easy to prove that from (\ref{tr23}),
(\ref{tr12}) and (\ref{tr24}) we get $\langle\G_{2},
\G_{4}\rangle= e, \langle\G_{1}, \G_{2}\rangle=e$ and
$\langle\G_{2}, \G_{3}\rangle=e$, respectively.

Relation (\ref{proj4a}) is now
\begin{eqnarray}
{}\ \ \ \ \ \ \ \ \G_{6'} \G_{5} \G_3 \G_4^{-1} \G_2 \G_5^{-2} \G_{4} \G_5^{2} \G_2^{-1} \G_4 \G_5^{-1} \G_4 \G_5^2 \G_{3}
\G_{2} \G_5^{-2} \G_4 \G_5^2  \G_2^{-1} \G_4 \G_{2}\G_3^{-2} \G_1 \G_3^{2} \G_1 =e. \label{project}
\end{eqnarray}

Now, equating the two expressions of $\G_{2'}$ given by (\ref{t4})
and (\ref{t5}), we get the relation
\begin{eqnarray}
\G_1^{-1} \G_2 \G_3^{-2} \G_{1} \G_3^{2} \G_2^{-1} \G_1 =
\G_{4}^{-1} \G_2 \G_5^{-2} \G_4  \G_{5}^2 \G_2^{-1}
\G_4,\label{kimat4}
\end{eqnarray}
which will be redundant later on.

\medskip

The relations which we have now are (\ref{66'}), (\ref{t3a}), (\ref{project}), (\ref{kimat4}) and
\begin{eqnarray}
{}\langle\G_{i}, \G_{j}\rangle & = & e, \ \ \mbox{$\G_{i}$ and $\G_{j}$ share a common triangle}\label{anibas1}\\
{}[ \G_{i} , \G_{j}] & = & e,  \ \ \mbox{$\G_{i}$ and $\G_{j}$ share no common triangle}\\
{}[ \G_{1} , \G_{3}^{-1} \G_{2} \G_{3}] = [ \G_{1} , \G_{3} \G_{2} \G_{3}^{-1}] & = & e \label{anibas3}\\
{}[ \G_{4} , \G_{5}^{-1} \G_{2} \G_{5}] = [ \G_{4} , \G_{5} \G_{2} \G_{5}^{-1}] &=& e.\label{anibas2}
\end{eqnarray}

Using (\ref{t3a}), we omit $\G_{6'}$, and therefore the group
$\pi_1(\C\P^2 \setminus B_4)$ has the generators
$\set{\G_i}_{i=1}^5$ and admits the relations (\ref{66'}),
(\ref{kimat4}), (\ref{anibas1}) - (\ref{anibas2}) for $i, j \neq
6'$, and the new form which  (\ref{project}) gets:
\begin{eqnarray}
{}\ \ \ \ \ \ \ \ \G_{5} \G_3 \G_2 \G_4^{-1} \G_2 \G_5^{-2} \G_{4} \G_5^{2} \G_2^{-1} \G_4 \G_5^{-1} \G_4 \G_5^2 \G_{3}
\G_{2} \G_5^{-2} \G_4 \G_5^2  \G_2^{-1} \G_4 \G_{2}\G_3^{-2} \G_1 \G_3^{2} \G_1  = e. \label{anilo}
\end{eqnarray}

\bigskip

Now we show that $\pi_1(\C\P^2 \setminus B_4)$  is isomorphic to a
quotient of $\B_6/\langle[X,Y]\rangle$, where $X, Y$ are
transversal half-twists. We choose a point in each triangle in
Figure \ref{home}. Then we choose a path $h_i$, connecting two
points in neighboring triangles, skipping the one which crosses
the edge $6$. We get a tree, see Figure \ref{tree}.
\begin{figure}[ht]
\begin{minipage}{\textwidth}
\begin{center}
\epsfbox{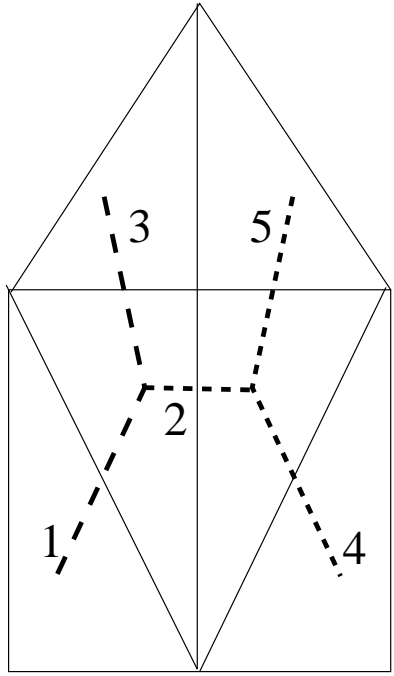}
\end{center}
\end{minipage}
\caption{The tree with five generators}\label{tree}
\end{figure}
The paths represent generators $\set{H_i}_{i=1}^5$ of the braid group
$\B_6$ with the following complete list of relations
\begin{eqnarray}
{}\langle H_{i}, H_{j}\rangle & = & e, \ \ \mbox{$H_{i}$ and $H_{j}$ are consecutive}\\
{}[ H_{i} , H_{j}] & = & e,  \ \ \mbox{$H_{i}$ and $H_{j}$ are disjoint}\\
{}[ H_{4} , H_{5} H_2 H_5^{-1}] & = & e \\
{}[ H_{1} , H_{3}^{-1} H_2 H_3] & = & e.
\end{eqnarray}

Denote $H_{6'} = {H}_{5} {H}_3 {H}_2 {H}_{3}^{-1} {H}_5^{-1}$,
where ${H}_{6'}$ corresponds to the missing path $h_6$,  being
transversal to $H_1$ and $H_2$ and disjoint from $H_4$.  Recall
the definition (\cite[Section~IV]{MTfg}) of the group
$\tilde{\B}_6$ as $\B_6/\langle[X,Y]\rangle$, and $X, Y$ are
transversal. Denote the images of $H_i$ as $\tilde{H}_i$ in
$\tilde{\B}_6$. Then the group $\tilde{\B}_6$ is generated by
$\tilde{H}_i$, $i=1, \dots, 5, 6'$, and the only relations are
\begin{eqnarray}
{}\langle\tilde{H}_{i}, \tilde{H}_{j}\rangle & = & e, \ \ \mbox{$\tilde{H}_{i}$ and $\tilde{H}_{j}$
are consecutive, $i, j \neq 6'$}\label{kimat6}\\
{}[\tilde{H}_{i} , \tilde{H}_{j}] & = & e,  \ \ \mbox{$\tilde{H}_{i}$ and $\tilde{H}_{j}$
are disjoint, $i, j \neq 6'$}\\
{}\ \ \ \ [\tilde{H}_{4} , \tilde{H}_{5} \tilde{H}_2
\tilde{H}_5^{-1}]=[\tilde{H}_{4} , \tilde{H}_{5}^{-1}
\tilde{H}_2 \tilde{H}_5] & = & e \\
{}\ \ \ \ \ \ \ \ [\tilde{H}_{1} , \tilde{H}_{3}^{-1} \tilde{H}_2
\tilde{H}_3]=[\tilde{H}_{1} , \tilde{H}_{3}
\tilde{H}_2 \tilde{H}_3^{-1}]& = & e\\
{}\ \ \ \ \ \ \ \ \tilde{H}_{5} \tilde{H}_3 \tilde{H}_2
\tilde{H}_{3}^{-1} \tilde{H}_5^{-1} & = &
\tilde{H}_{6'},\label{kimat7}
\end{eqnarray}
where $\tilde{H}_{4}$ and $\tilde{H}_{5}^{-1} \tilde{H}_2
\tilde{H}_5$ ($\tilde{H}_{1}$ and  $\tilde{H}_{3} \tilde{H}_2
\tilde{H}_3^{-1}$, respectively) are transversal. We note that
(\ref{kimat7}) can be used to remove $\tilde{H}_{6'}$ from the
list of generators (in the same way as $\Gamma_{6'}$ has been
eliminated from the presentation of $\pi_1(\C\P^2 \setminus
B_4)$).

According to our result, $\pi_1(\C\P^2 \setminus B_4)$ is a
quotient of $\tilde{\B}_6$. Now we eliminate (\ref{kimat4}). Since
$\G_3^{-2} \G_{1} \G_3^{2}$ and $\G_3^{-1} \G_{2} \G_3$ are
transversal, the relations in   $\tilde{\B}_6$ imply that they
commute and so the left hand side of (\ref{kimat4}) is equal to
\begin{eqnarray*}
{}\ \ \ \G_1^{-1} \G_2 (\G_3^{-1} \G_2^{-1} \G_3) \G_3^{-2} \G_{1}
\G_3^2 (\G_3^{-1} \G_2 \G_{3}) \G_2^{-1} \G_1 =  \G_1^{-1}
\G_3^{-1} \G_2^{-1} \G_{1} \G_2 \G_3 \G_1 = \G_2.
\end{eqnarray*}

Similarly, the right hand side of (\ref{kimat4}) is also equal to
$\G_2$. This allows us to eliminate (\ref{kimat4}). Since both
sides of (\ref{kimat4}) are equal to $\G_2$, we have shown that
$\G_2=\G_{2'}$, therefore $\G_6=\G_{6'}$ (see (\ref{t3}) and
(\ref{t3a})). That means that (\ref{66'}) is redundant too.  Thus
$\pi_1(\C\P^2 \setminus~B_4)$ is isomorphic to
$\tilde{\B}_6/\langle(\ref{anilo})\rangle$.

\medskip

In order to get the group $\Pi_{(B_4)}$, we take $\G_j^2=e$ for each $j$. Relation (\ref{anilo})
is then redundant. By \cite{RTV},
the rest of the relations in $\pi_1(\C\P^2 \setminus B_4)$, together with the ones $\G_j^2=e$,
are the only ones which are required in order to make $\Pi_{(B_4)}$ isomorphic to $S_{6}$.
\end{proof}

\forget

{\bf Appendix.}\\
We show that $[ \G_{2}, \G_{6} ]=e$ from (\ref{com26a}) gets the
form $[\G_{4} , \G_{5}^{-1} \G_{2} \G_{5}]=e$ and relation
(\ref{66'}) is redundant:
\begin{eqnarray*}
{}& &e=[\G_{2}, \G_{6}]\stackrel{(\ref{t3})}{=}[\G_{2} , \G_{5}
\G_{3} \G_{2'} \G_{3}^{-1} \G_{5}^{-1}]= [\G_{5}^{-1} \G_{2}
\G_{5}, \G_{3} \G_{2'} \G_{3}^{-1}]\stackrel{(\ref{t4})}{=}\\ & &
[\G_{5}^{-1} \G_{2} \G_{5}, \G_{3} \G_{4}^{-1}\G_{2} \G_{4'}
\G_2^{-1} \G_4 \G_{3}^{-1}] \stackrel{(\ref{com346}),
(\ref{anibas2}), (\ref{25tr})}{=} [\G_{2}^{-1} \G_{3}^{-1} \G_{2}
\G_{5} \G_{2}^{-1}\G_{3} \G_{2},  \G_{4'}]
\stackrel{(\ref{tr24})}{=}\\ & & [\G_{3}\G_2^{-1} \G_{3}^{-1}
\G_{5} \G_3\G_2\G_{3}^{-1}, \G_{4'}] \stackrel{(\ref{com346}),
(\ref{35com})}{=}[\G_2^{-1}\G_{5}\G_2, \G_{4'}]
\stackrel{(\ref{t4})}{=}\\& &[\G_2^{-1}\G_{5}\G_2,\G_{2}^{-1}
\G_{2'}^{-1} \G_4 \G_{2'} \G_2]
\stackrel{(\ref{tr23})}{=}[\G_{5},\G_4 \G_{2'}\G_4^{-1}]
\stackrel{(\ref{t5})}{=}\\& &[\G_{5},\G_4 \G_{1}^{-1} \G_2 \G_{1'}
\G_2^{-1}\G_{1}\G_4^{-1}] \stackrel{(\ref{com14}), (\ref{com15}),
(\ref{tr12})}{=} [\G_{5},\G_4 \G_{1'}^{-1} \G_2
\G_{1'}\G_4^{-1}]=\\& &[\G_{5},\G_4  \G_2 \G_4^{-1}]
\stackrel{(\ref{tr23}), (\ref{25tr})}{=}[\G_{4},\G_5^{-1}  \G_2
\G_5],
\end{eqnarray*}

\smallskip

\begin{eqnarray*}
{}& &e=[\G_{6}, \G_{6'}]\stackrel{(\ref{t3}),
(\ref{t3a})}{=}[\G_{5} \G_{3} \G_{2'} \G_{3}^{-1} \G_{5}^{-1},
\G_{5} \G_{3} \G_{2} \G_{3}^{-1} \G_{5}^{-1}]=[\G_2,
\G_{2'}]\stackrel{(\ref{t4})}{=}\\ & & [\G_{2}, \G_{4}^{-1}\G_{2}
\G_{4'} \G_2^{-1} \G_4]= [\G_{2}^{-1}\G_{4}\G_{2} \G_{4}^{-1}
\G_2, \G_{4'}] \stackrel{(\ref{t2})}{=}[\G_{4}\G_{2} \G_{4}^{-1}
\G_{2}\G_{4}\G_{2}^{-1} \G_{4}^{-1}, \G_{5}^{-1}\G_{5}^{-1}\G_{4}
\G_{5} \G_5]\stackrel{(\ref{anibas2}), (\ref{tr35})}{=}\\ & &
[\G_5\G_{2}\G_{5}^{-1},\G_{4}\G_{2}^{-1} \G_{4}^{-1} \G_4
\G_{5}\G_4^{-1}\G_{4}\G_{2} \G_{4}^{-1} ]
\stackrel{(\ref{anibas2})}{=}[\G_5\G_{2}\G_{5}^{-1},\G_{2}^{-1}\G_{5}\G_{2}]
\stackrel{(\ref{25tr})}{=}\\ & &
[\G_5\G_{2}\G_{5}^{-1},\G_{5}\G_{2}\G_5^{-1}].
\end{eqnarray*}

\forgotten


\begin{thebibliography}{99}


\bibitem{A}
Amram, M., {\em Galois Covers of Algebraic Surfaces,} Ph.D. dissertation, 2001.



\bibitem{AT-TT1}
Amram, M., Teicher, M., {\em On the degeneration, regeneration and braid monodromy
of the surface $T \times T$}, Acta Applicandae Mathematicae, 75(1) (2003), 195-270.



\bibitem{AT-TT2}
Amram, M., Teicher, M., {\em The fundamental group of  the complement of the
branch curve of the surface $T \times T$}, OJM, Japan, 40(4) (2003), 1-37.



\bibitem{K3-a}
Amram, M., Ciliberto, C., Miranda, R., Teicher, M.,
{\em Braid monodromy factorization for a non-prime $K3$ surface branch curve}, submitted.




\bibitem{Gol1}
Amram, M., Goldberg, D., Teicher, M., Vishne U., {\em The fundamental group of a
Galois cover of $\C\P^1 \times T$}, Algebraic and Geometric Topology 2 (2002), 403-432.



\bibitem{Gol2}
Amram, M., Goldberg, D., {\em Higher degree Galois covers of $\C\P^1 \times T$},
Algebraic and Geometric Topology 4 (2004), 841-859.


\bibitem{Cox}
Amram, M., Teicher, M., Vishne, U., {\em The Coxeter quotient of the
fundamental group of a Galois cover of $T \times T$}, Communications in Algebra 34 (2006), 89-106.



\bibitem{ATV}
Amram, M., Teicher, M., Vishne, U., {\em The fundamental group of the Galois
cover of the surface $T \times T$}, submitted.



\bibitem{ATV-H}
Amram, M., Teicher, M., Vishne, U., {\em Fundamental groups of Galois covers
of Hirzebruch surfaces $F_1(2,2)$}, submitted.




\bibitem{artin}
Artin, E., {\em Theorie der Z\"opfe}, Abhandlungen aus dem Mathematischen Seminar
der Hamburgischen Universit\"at, vol. 4, 1926.




\bibitem{denis}
Auroux, D., Donaldson, S. K., Katzarkov, L., Yotov, M.,
{\em Fundamental groups of complements of plane curves and symplectic invariants},
Topology 43 (2004), 1285-1318.



\bibitem{AK}
Auroux, D., Katzarkov, L., {\em Branched coverings of $\C\P^2$ and invariants of symplectic
$4$-manifolds}, Invent. Math. 142 (2000), 631-673.




\bibitem{BGT}
Bruns, W., Gubeladze, J., Trung, N.V.,  \newblock {\em Normal polytopes,
triangulations, and Koszul algebras}
\newblock J. reine angew. Math. 485 (1997), 123-160.




\bibitem{F}
Fulton, W.,  \newblock {\em Introduction to Toric Varieties}, Ann. of Math.
Studies, vol. 131,  Princeton Univ. Press, 1993.



\bibitem{Hu}
Hu, S.,  \newblock {\em Semi-stable degeneration of toric varieties and their
hypersurfaces}, 2003, AG/0110091.



\bibitem{mag}
Magnus, W.,  {\em \"Uber Automorphismen von Fundamentalgruppen
berandeter Fl\"achen}, Mathematische Annalen 109 (1934), 617-646.





\bibitem{MII}
Moishezon, B., {\em Algebraic surfaces and the
arithmetics of braids, II}, Cont. Math. 44 (1985), 311-344.




\bibitem{M2} Moishezon, B., {\em On cuspidal branch curves}, J. Algebraic Geometry 2 (1993), 309-384.




\bibitem{MRT}
Moishezon, B., Robb, A.,  Teicher, M., {\em On Galois covers of Hirzebruch
surfaces}, Math. Ann. 305 (1996), 493-539.


\bibitem{MT1}
Moishezon, B.,  Teicher, M., {\em  Simply connected algebraic surfaces of positive
index}, Invent. Math. 89  (1987), 601-643.


\bibitem{MT2}
 Moishezon, B., Teicher, M., {\it Braid group technique in complex geometry I,
Line arrangements in $\C\P^2$}, Contemporary Math. 78 (1988),  425-555.


\bibitem{MT3}
Moishezon, B., Teicher, M., {\it Braid group
   technique in complex geometry II, From arrangements of lines and conics
   to cuspidal curves}, Algebraic Geometry, Lect. Notes in Math., vol. 1479, 1990, 131-180.



\bibitem{ZcZ}
Moishezon, B., Teicher, M., {\it Finite fundamental
groups, free over $\Z/c\Z$, Galois covers of $\C \P^2$},
Math. Ann. 293 (1992), 749-766.



\bibitem{MT4}
Moishezon, B., Teicher, M., {\it  Braid group
technique in complex geometry III: Projective degeneration of $V_3$, }
Contemp. Math. 162 (1993), 313-332.



\bibitem{MT5}
Moishezon, B.,  Teicher, M., {\it Braid group technique in complex geometry IV:
Braid monodromy of the branch curve $S_3$ of $V_3 \rightarrow \C\P^2$ and application
to $\pi_1(\C\P^2 - S_3, \ast)$,} Contemporary Math. 162 (1993),  332-358.


\bibitem{MTfg}
Moishezon, B.,  Teicher, M., {\it Braid group technique in complex geometry V:
The fundamental group of a complement of a branch curve of a Veronese generic projection},
Communications in Analysis and Geometry 4 (1996), 1-120.



\bibitem{Od}
T.~Oda,    \newblock {\em Convex Bodies and Algebraic Geometry,}  Springer-Verlag, 1985.



\bibitem{robb}
Robb, A., {\it The topology of branch curves of complete intersections}, Doctoral Thesis,
Columbia, 1994.



\bibitem{RTV}
Rowen, L.H., Teicher, M., Vishne, U., {\it Coxeter Covers of the Symmetric Groups},
J. Group Theory 8 (2005), 139-169.


\bibitem{Str}
B.~Sturmfels,  \newblock {\em Equations defining toric varieties,}
Algebraic Geometry-Santa Cruz 1995, Proc. Sympos. Pure Math. 62 (1997), 437-449.


\bibitem{vK}
\vK, E.R., {\em On the fundamental group of an algebraic curve},
Amer. J. Math. 55 (1933), 255-260.



\bibitem{Z2}
Zariski, O., {\em On the problem of existence of algebraic functions of two
variables possesing a given branch curve}, American Journal of Math., vol. 51, 1929.



\bibitem{Z}
Zariski, O., {\em On the Poincar\'e group of rational plane curve},
Amer. J. Math. 58 (1936), 607-619.





\end{thebibliography}
\end{document}